	\magnification=\magstep1  
	\overfullrule=0pt  
	\def\lb{\lbrack}	\def\rb{\rbrack}  
	\def\IN{\mathop{{\rm I}\kern-.2em{\rm N}}\nolimits}  
	\def\IR{\mathop{{\rm I}\kern-.2em{\rm R}}\nolimits}  
	\def\A{{\cal A}}  
	\def\C{{\cal C}}  
	\def\D{{\cal D}}  
	\def\E{{\cal E}}  
	\def\F{{\cal F}}  
	\def\bone{{\bf 1}}  
	\def\osc{\mathop{\rm osc}\nolimits}  
	\def\diam{\mathop{\rm diam}\nolimits}  
	\def\rf#1{{\rm $\lbrack$#1$\rbrack$}}  
	\def\iitem{\itemitem}  
	\def\ov{\overline}  
	\def\varep{\varepsilon}  
	\def\frac#1#2{{\textstyle{#1\over#2}}}  
	\def\subsetneqq{\mathop{\ \lower.5ex\hbox{$\buildrel{\scriptstyle   
		\subset}\over{\scriptstyle\ne}$}\ }\nolimits}  
	\def\supsetneqq{\mathop{\ \lower.5ex\hbox{$\buildrel{\scriptstyle   
		\supset}\over{\scriptstyle\ne}$}\ }\nolimits}  
	\def\trivert{\mathop{\ \hbox{$|\kern-.4em|\kern-.4em|$}\ }\nolimits}  
	\def\blackbox{\hbox{\vrule width6pt height7pt depth1pt}}  
	\def\qed{\hfill~\blackbox}  
	\def\refrule{\vrule height0.4pt depth0pt width1in}  
	\def\proof{\noindent{\it Proof.}\quad}  
	\outer\def\demo#1. #2\par{\medbreak\noindent {\it#1.\enspace}  
		{\rm#2}\par\ifdim\lastskip<\medskipamount\removelastskip  
		\penalty55\medskip\fi}  
\topinsert\vskip.5truein\endinsert  
\centerline{\bf On certain classes of Baire-1 functions}  
\smallskip  
\centerline{\bf with applications to Banach space theory}  
\bigskip  
\centerline{R. Haydon, E. Odell and  H. Rosenthal  
	\footnote{}{Research partially supported by NSF Grant DMS-8601752.}}  
\vskip.3in  
\centerline{\bf Abstract}  
\smallskip  
Certain subclasses of $B_1(K)$, the Baire-1 functions on a compact  
metric space $K$, are defined and characterized. Some applications 
to Banach spaces are given. 
\vskip.3in 

\baselineskip=18pt		  

\beginsection{0. Introduction.}  

Let $X$ be a separable infinite dimensional Banach space and let $K$ denote  
its dual ball, $Ba (X^*)$, with the weak* topology. $K$ is compact metric  
and $X$ may be naturally identified with a closed subspace of $C(K)$.   
$X^{**}$ may also be identified with a closed subspace of $A_\infty (K)$,  
the Banach space of bounded affine functions on $K$ in the sup norm.  
Our general objective is to deduce information about the isomorphic   
structure of $X$ or its subspaces from the topological nature of the functions  
$F\in X^{**} \subseteq A_\infty (K)$. A classical example of this type  
of result is: $X$ is reflexive if and only if $X^{**} \subset C(K)$.  

A second example is the following theorem. ($B_1 (K)$ is the class of  
bounded Baire-1 functions on $K$ and $DBSC (K)$ is the subclass of differences  
of bounded semicontinuous functions on $K$. The precise definitions  
appear below in \S1.) We write $Y\hookrightarrow X$ if $Y$ is isomorphic  
to a subspace of $X$.  

\proclaim Theorem A.  
Let $X$ be a separable Banach space and let $K= Ba (X^*)$ with the 
weak* topology.  
\smallskip  
\iitem{a)} \rf{35} $\ell_1 \hookrightarrow X$ iff $X^{**}\setminus  
B_1 (K) \ne \emptyset$.  
\iitem{b)} \rf{7} $c_0 \hookrightarrow X$ iff $\lbrack X^{**}\,\cap\,  
DBSC (K)\rbrack \setminus C(K)\ne \emptyset$.  
\smallskip\par  

Theorem A provides the motivation for this paper: What can be said   
about $X$ if $X^{**} \,\cap\, \lbrack B_1 (K)\setminus DBSC(K)\rbrack  
\ne \emptyset$?  
To study this problem we consider various subclasses of $B_1(K)$  
for an arbitrary compact metric space $K$. J.~Bourgain   
has also used this approach and some of our results and techniques  
overlap with those of \rf{8,9,10}. In a different direction,   
generalizations of $B_1(K)$ to spaces where $K$ is not compact metric  
with ensuing applications to Banach space theory have been developed  
in \rf{22}.   

In \S1 we consider two subclasses of $B_1(K)$ denoted $B_{1/4}(K)$  
and $B_{1/2} (K)$ satisfying  
$$C(K) \subseteq DBSC (K) \subseteq B_{1/4} (K) \subseteq B_{1/2}(K)  
\subseteq B_1(K)\ .\leqno(0.1)$$  
Our interest in these  classes stems from Theorem B (which we  
prove in \S3).  

\proclaim Theorem B.  
Let $K$ be a compact metric space and let $(f_n)$ be a uniformly bounded  
sequence in $C(K)$ which converges pointwise to $F\in B_1 (K)$.  
\smallskip  
\iitem{a)} If $F\notin B_{1/2}(K)$, then $(f_n)$ has a subsequence whose  
spreading model is equivalent to the unit vector basis of $\ell_1$.  
\iitem{b)} If $F\in B_{1/4} (K) \setminus C(K)$, there exists  
$(g_n)$, a convex block subsequence of $(f_n)$, whose spreading model  
is equivalent to the summing basis for $c_0$.   
\smallskip\par  

Theorem B may be regarded as a local version of Theorem~A  
(see Corollary~3.10). In fact the proof is really a localization of  
the proof of Theorem~A.  In Theorem~3.7 we show that the converse to  
a) holds and thus we obtain a characterization of $B_1(K)\setminus  
B_{1/2} (K)$ in terms of $\ell_1$ spreading models. We do not know  
if the condition in b) characterizes $B_{1/4} (K)$ (see Problem~8.1).  

Given that our main objective is to deduce information about the subspaces  
of $X$ from the nature of $F\in X^{**} \,\cap\, B_1 (K)$, it is   
useful to introduce the following definition.  

Let $\C$ be a class of separable   
infinite-dimensional Banach spaces and let $F\in B_1 (K)$.  
$F$ is said to {\it govern\/} $\C$ if whenever $(f_n) \subseteq C(K)$  
is a uniformly bounded sequence converging pointwise to $F$, then  
there exists a $Y\in \C  $ which embeds into   
$[(f_n)]$, the closed linear span of $(f_n)$.  
We also say that $F$ {\it strictly governs\/} $\C$ if whenever $(f_n)  
\subseteq C(K)$ is a uniformly bounded sequence converging pointwise  
to $F$, there exists a convex block subsequence $(g_n)$ of $(f_n)$  
and a $Y\in \C$ with $\lbrack (g_n)\rbrack $ isomorphic to $Y$.  

Theorem A (b) can be more precisely formulated as: if $F\in DBSC (K)  
\setminus C(K)$, then $F$ governs $\{ c_0\}$.  (In fact Corollary~3.5  
below yields that $F \in B_1(K) \setminus C(K)$ strictly governs  
$\{c_0\}$ if and only if $F\in DBSC(K)$.)  
In \S4 we prove that  
the same result holds if ${F\in DSC (K)\setminus C(K)}$. (A more  
general result, with a different proof, has been obtained by Elton \rf{13}.)  
We also note in \S4 that there are functions that govern $\{ c_0\}$  
but are not in $DSC(K)$.   

In \S6 we give a characterization of $B_{1/4}(K)$ (Theorem~6.1) and use  
it to give an example of an $F\in B_{1/4} (K)\setminus C(K)$ which  
does not govern $\{c_0\}$. Thus Theorem~B (b) is best possible.  

In \S7 we note that there exists a $K$ and an $F\in B_{1/2}(K)$ which governs  
$\{ \ell_1\}$. We also give an example of an $F\in B_{1/2}(K)$ which  
governs $\C = \{ X: X$ is separable and $X^*$ is nonseparable$\}$ but does  
not govern $\{ \ell_1\}$.  

\S1 contains the definitions of the classes $DBSC(K)$, $DSC(K)$,  
$B_{1/2}(K)$ and $B_{1/4}(K)$. At the end of \S1 we briefly recall  
the notion of spreading model. In \S2 we recall some ordinal indices  
which are used to study $B_1(K)$. A detailed study of such indices  
can be found in \rf{25}. Our use of these indices and many of the  
results of this paper have been motivated by \rf{8,9,10}.  
Proposition~2.3 precisely characterizes $B_{1/2}(K)$ in terms of our index.  

In \S5 we show that the inclusions in (0.1) are, in general, proper.  
We first deduce this from a Banach space perspective. Subsequently, we  
consider the case where $K$ is countable.  Proposition~5.3 specifies  
precisely  
how large $K$ must be in order for each separate inclusion in (0.1)  
to be proper.  

In \S8 we summarize some problems raised throughout this paper and raise  
some new questions regarding $B_{1/4}(K)$.  

We are hopeful that our approach will shed some light on the central problem:  
if $X$ is infinite dimensional, does $X$ contain an infinite dimensional  
reflexive subspace or an isomorph of $c_0$ or $\ell_1$?  A different  
attack has been mounted on this problem in the last few years by   
Ghoussoub and Maurey. The interested reader should also consult their   
papers ({\it e.g.}, \rf{18,19,20,21}).  Another fruitful approach has  
been via the theory of types (\rf{26}, \rf{24}, \rf{38}). We wish to  
thank S.~Dilworth and R.~Neidinger for useful suggestions.  

\beginsection{1. Definitions.}  

In this section we give the basic definitions of the Baire-1 subclasses  
in which we are interested. Let $K$ be a compact metric space. 
$B_1(K)$ shall denote the class of bounded Baire-1 functions on $K$,  
{\it i.e.}, the pointwise limits of (uniformly bounded) pointwise  
converging sequences $(f_n) \subseteq C(K)$. $DBSC(K) = \{ F: K\to \IR\mid$  
there exists $(f_n)_{n=0}^\infty \subseteq C(K)$ and $C<\infty$  
such that $f_0 \equiv 0$, $(f_n)$ converges pointwise to $F$ and  
$$\sum_{n=0}^\infty |f_{n+1}(k) - f_n(k)| \le C\ \hbox{ for all }\   
k\in K\bigr\}\ .\leqno(1.1)$$  
If $F\in DBSC(K)$ we set $|F|_D = \inf \{ C\mid $ there exists  
$(f_n)_{n=0}^\infty \subseteq C(K)$ converging pointwise to $F$   
satisfying (1.1)  
with $f_0 \equiv 0\}$.  $DBSC(K)$ is thus precisely  
those $F$'s which are the ``difference of bounded semicontinuous  
functions on $K$.'' Indeed if $(f_n)$ satisfies (1.1), then  
$F= F_1 -F_2$ where $F_1(k) = \sum_{n=0}^\infty (f_{n+1} -f_n)^+ (k)$  
and $F_2(k) = \sum_{n=0}^\infty (f_{n+1}-f_n)^- (k)$ are both (lower)  
semicontinuous. The converse is equally trivial.  

It is easy  to prove that ($DBSC(K), |\cdot|_D$) is a Banach space by using  
the series criterion for completeness. The fact that $\Vert F\Vert_\infty  
\le |F|_D$ but the two norms are in general not equivalent on $DBSC(K)$, leads  
naturally to the following two definitions.  
$$\eqalign{ B_{1/2} (K) &= \bigl\{ F\in B_1 (K)\mid \hbox{ there exists a  
	sequence}\cr  
&\qquad (F_n) \subseteq DBSC(K)\ \hbox{ converging uniformly to $F$}\bigr\}  
	\hbox{ and}\cr  
B_{1/4}(K) & = \bigl\{ F\in B_1 (K) \mid \hbox{ there exists } (F_n)\cr  
&\qquad \hbox{converging uniformly to $F$ with } \sup_n |F_n|_D <\infty  
\bigr\}\ .\cr}$$  

It can be shown that $DBSC(K)$ is a Banach algebra  under pointwise  
multiplication, and hence $B_{1/2} (K)$ can be identified with  
$C(\Omega)$, where $\Omega$ is the ``structure space'' or ``maximal  
ideal space'' of $\Omega$. Thus $B_{1/4}(K)$ also has a natural   
interpretation in the general context of commutative Banach algebras.  

There is a natural norm on $B_{1/4}(K)$ given by  
$$|F|_{1/4} = \inf \bigl\{ C :  \hbox{ there exists $(F_n)$ converging   
	uniformly with } \sup_n |F_n|_D \le C\bigr\}$$  
Furthermore $(B_{1/4}(K),|\cdot|_{1/4})$ is a Banach space.  One way  
to see this is to use the following elementary  

\proclaim Lemma 1.1.  
Let $(M,d_1)$ be a complete metric space and let $d_2$ be a metric on $M$  
with $d_1 (x,y) \le d_2(x,y)$ for all $x,y\in M$. If all $d_2$-closed  
balls in $M$ are also $d_1$-closed, then $(M,d_2)$ is complete.\par  

The hypotheses of the lemma apply to $M= \{ F:|F|_{1/4} \le 1\}$ and  
$d_1,d_2$ given, respectively, by $\Vert\cdot\Vert_\infty$    
and $|\cdot|_{1/4}$.  

\demo Remark 1.2.  
While we shall confine our attention to $B_{1/2}$ and $B_{1/4}$, one could   
of course continue the game, defining   
$$\eqalign{ B_{1/8} (K) & = \bigl\{ F\in B_1(K)  \mid \hbox{ there exists }  
	(F_n)\subseteq DBSC(K) \cr  
&\qquad \hbox{with } \ |F_n-F|_{1/4} \to 0\bigr\}\ \hbox{ and }\cr  
B_{1/16} (K) & = \bigl\{ F\in B_1(K) \mid \hbox{ there exists } F_n\cr  
&\qquad \hbox{with }\ \sup_n |F_n|_D <\infty\ \hbox{ and }\   
	|F_n-F|_{1/4} \to 0\bigr\}\ .\cr}$$  
This could be continued obtaining  
$$DBSC(K) \subseteq \cdots \subseteq B_{1/2^{2n}}(K) \subseteq  
B_{1/2^{2n-1}}(K) \subseteq \cdots \subseteq B_{1/2} (K)$$  
with $B_{1/2^{2n}}(K)$ having a norm $|\cdot|_{1/2^{2n}}$ which,   
using Lemma~1.1, is easily seen to be complete.  

There is another class of Baire-1 functions that shall interest us, the  
differences of (not necessarily bounded) semi-continuous functions  on $K$.  
$$\eqalign{DSC(K) & = \bigl\{ F: K\to \IR\mid \hbox{ there exists a uniformly  
	bounded sequence}\cr  
&\qquad (f_n)_{n=0}^\infty \subseteq C(K)\ \hbox{ converging pointwise to $F$  
	with}\cr  
&\qquad \sum_{n=0}^\infty |f_{n+1} (k) - f(k)| < \infty \ \hbox{ for }\   
	k\in K\bigr\}\ .\cr}$$  
An interesting subclass of $DSC(K)$, is $PS(K)$, the pointwise limits  
of {\it pointwise stabilizing\/} (pointwise ultimately constant) sequences.  
$$\eqalign{PS(K) & = \bigl\{ F\in B_1 (K)\mid \hbox{ there exists a   
	uniformly bounded sequence}\cr  
&\qquad (f_n)\subseteq C(K) \ \hbox{ with the property that for all $k\in K$  
	there exists}\cr  
&\qquad m\in \IN\ \hbox{ such that }\ f_n(k) = F(k)\ \hbox{ for }\   
	n\ge m\bigr\}\ .\cr}$$  

\demo Remark 1.3.  
We discuss $PS(K)$ in Proposition~4.9. Both of these classes were considered  
in \rf{10}, and as noted there, if an indicator function    
$\bone_A \in B_1 (K)$, then $\bone_A \in PS (K)$. Indeed $A$ must be both  
$F_\sigma$ and $G_\delta$ (cf.\ Proposition 2.1 below)  
and so we can write $A= \bigcup_n F_n =   
\bigcap_n G_n$ where $F_1\subseteq F_2\subseteq \cdots\ $ are closed sets  
and $G_1 \supseteq G_2\supseteq \cdots\ $ are open sets. Then by the  
Tietze extension theorem, for each $n$ choose $f_n\in Ba (C(K))$ with  
$f_n$ identically 1 on $ F_n$ and identically $0$ on $K\setminus G_n$.  
Thus for all $k\in K$, $(f_n(k))_n$ is ultimately $\bone_A (k)$.  

The {\it summing basis\/} $(s_n)$ for (an isomorph of) $c_0$ is  
characterized by   
$$\Vert \sum a_ns_n\Vert = \sup_k |\sum_{i=1}^k a_i|\ .$$  

Let $(x_n)$ be a seminormalized basic sequence. A basic sequence  
$(e_n)$ is said to be a {\it spreading model\/} of $(x_n)$  
if for all $k\in \IN$ and all $\varep >0$ there exist $N$ so that  
if $N<n_1 < n_2 < \cdots  < n_k$ and $(a_i)_1^k \subseteq \IR$  
with $\sup_i |a_i| \le 1$, then  
$$\Big| \Vert \sum_{i=1}^k a_i x_{n_i} \Vert -  
\Vert \sum_{i=1}^k a_ie_i\Vert \Big| <\varep\ .$$  
For further information on spreading models see  
\rf{4}.  

We recall that if $(f_n) \subseteq Ba(C(K))$ converges pointwise to  
$F\in B_1 (K)\setminus C(K)$ then there exists a $C= C(F)$ such that  
$(f_n)$ has a basic subsequence $(f'_n)$ with basis constant $C$  
which $C$-dominates $(s_n)$. Thus $C\Vert \sum a_nf'_n\Vert \ge \Vert  
\sum a_ns_n\Vert$, for all $(a_n) \subseteq \IR$ (see {\it e.g.},  
\rf{31}). Furthermore $(f'_n)$ can be taken to have a spreading model  
\rf{4}. The constant $C$ depends only on $\sup \{ \osc (F,k) \mid   
k\in K\}$ (see \S2 for the definition of $\osc (F,k)$).  

Finally we recall that a sequence $(g_n)$ in a Banach space is a  
{\it convex block subsequence\/} of $(f_n)$ if $g_n = \sum_{i=p_n+1}^  
{p_{n+1}} a_if_i$ where $(p_n)$ is an increasing sequence of integers,  
$(a_i) \subseteq \IR^+$ and for each $n$, $\sum_{i= p_n+1}^{p_{n+1}}  
a_i = 1$.  

\beginsection{2. Ordinal Indices for $B_1 (K)$.}  

Let $(K,d)$ be a compact metric space and let $F:K\to \IR$ be a bounded  
function. The Baire characterization theorem \rf{3} states that  
$F\in B_1(K)$ iff for all closed nonempty $L\subseteq K$, $F\big|_L$ has a   
point of continuity (relative to the compact space $(L,d)$). This leads  
naturally to an ordinal index for Baire-1 functions which we now  
describe.  

For a closed set $L\subseteq K$ and $\ell \in L$   
let the {\it oscillation of\/}  
$F\big|_L$ at $\ell$ be given by 
$\osc_L (F,\ell)= \lim_{\varep\downarrow 0} \sup  
\{ f(\ell_1) - f(\ell_2) \mid \ell_i\in L$ and $d(\ell_i,\ell) <\varep$  
for $i= 1,2\}$. We define  the {\it oscillation of $F$ over\/} $L$ by  
$\osc_L F = \sup \{ F(\ell_1) - F(\ell_2) \mid \ell_1,\ell_2\in L\}$.  

For $\delta >0$, let $K_0 (F,\delta) = K$ and if $\alpha < \omega_1$ let  
$$K_{\alpha +1} (F,\delta) = \bigl\{ k\in K_\alpha (F,\delta) \mid  
\osc_{K_\alpha (F,\delta)} (F,k) \ge \delta\bigr\}\ .$$  
For limit ordinals $\alpha$, set  
$$K_\alpha (F,\delta) = \bigcap_{\beta<\alpha} K_\beta (F,\delta)\ .$$  
Note that $K_\alpha (F,\delta)$ is always closed and $K_\alpha (F,\delta)  
\supseteq K_\beta (F,\delta)$ if $\alpha <\beta$.  The index $\beta (F,  
\delta)$ is given by  
$$\beta (F,\delta) = \inf \bigl\{   \alpha <\omega_1 \mid K_\alpha (F,  
\delta) = \emptyset\bigr\}$$  
provided $K_\alpha (F,\delta) = \emptyset$ for some $\alpha <\omega_1$  
and $\beta (F,\delta) = \omega_1$ otherwise. Since $K$ is separable,  
the transfinite sequence $(K_\alpha (F,\delta))_{\alpha <\omega_1}$  
must stabilize: there exists $\beta < \omega_1$ so that   
$K_\alpha (F,\delta) = K_\beta (F,\delta)$ for $\beta \ge \alpha$.  

The Baire characterization theorem yields that $\beta (F,\delta) <  
\omega_1$ for all $\delta >0$ iff $F\in B_1(K)$. In fact we have the  
following proposition. In its statement $\A$ denotes the  algebra  
of {\it ambiguous\/} subsets of $K$. Thus $A\in \A$ iff $A$ is both  
$F_\sigma$ and $G_\delta$. Also we write $\lb F\le a\rb$ for the set  
$\{ k\in K\mid F(k)\le a\}$.  

\proclaim Proposition 2.1.  
Let $F: K\to \IR$ be a bounded function on the compact metric   
space $K$. The following are equivalent.  
\smallskip  
\iitem{1)} $F\in B_1 (K)$.  
\iitem{2)} $\beta (F,\delta) < \omega_1$ for all $\delta >0$.  
\iitem{3)} For $a$ and $b$ real, $\lb F\le a\rb$   
and $\lb F\ge b\rb$ are both   
$G_\delta$ subsets of $K$. 
\iitem{4)} For $U$ an open subset of $\IR$, $F^{-1}(U)$ is an 
$F_\sigma$ subset of $K$.   
\iitem{5)} For $a<b$, $\lbrack F\le a\rbrack$     
and $\lbrack F\ge b\rbrack$ may be separated by  
disjoint sets in $\A$. Equivalently, there exists $A\in \A$ with   
$\lbrack F\le a\rbrack \subseteq A$ and $A\,\cap\, \lbrack F\ge b\rbrack     
=\emptyset$.  
\iitem{6)} $F$ is the uniform limit of a sequence of $\A$-simple  
functions ($\A$-measurable functions with finite range).  
\iitem{7)} $F$ is the uniform limit of a sequence $(g_n) \subseteq DSC(K)$.  
\iitem{8)} $F$ is the uniform limit of a sequence $(g_n)\subseteq PS(K)$.  
\smallskip\par  

The proof is standard and can be compiled from \rf{23}. We are more interested  
in an analogous characterization of $B_{1/2} (K)$. Before stating that  
proposition we need a few more definitions.  

$\D$ shall denote the algebra of all finite unions of differences of   
closed subsets of $K$. $\D$ is easily seen to be a subalgebra of $\A$.  

One of the statements in our next proposition involves another  
ordinal index for Baire-1 functions, $\alpha (F;a,b)$, which as we  
shall see is closely related to our index. For $a<b$, let $K_0 (F;a,b)=K$  
and for any ordinal $\alpha$, let   
$$\eqalign{K_{\alpha +1} (F;a,b) &= \{ k\in K_\alpha  
(F;a,b) \mid \ \hbox{ for all }\ \varep >0 \hbox{ and } i=1,2,\cr  
&\qquad \hbox{there exist } k_i \in K_\alpha  
	(F;a,b) \ \hbox{ with }\ d(k_i,k)\le\varep\ ,\cr  
&\qquad  F(k_1) \ge b \ \hbox{ and }\ F(k_2) \le a\}\ .\cr}$$  
Equivalently, $K_{\alpha +1} = \ov{K_\alpha\,\cap\,     
\lbrack F\le a\rbrack }\,\cap\,  
\ov{K_\alpha\,\cap\, \lbrack F\ge b\rbrack}$.     
At limit ordinals $\alpha$ we set  
$$K_\alpha (F;a,b) = \bigcap_{\beta <\alpha} K_\beta (F;a,b)\ .$$  
As before these sets are closed and decreasing.  We let $\alpha (F;a,b)  
=\inf \{ \gamma <\omega_1 \mid K_\gamma (F;a,b) =\emptyset\}$ if  
$K_\gamma (F;a,b) =\emptyset$ for some $\gamma <\omega_1$ and  
let $\alpha (F;a,b) =\omega_1$ otherwise.   

\demo Remark 2.2.  
The index $\alpha (F;a,b)$ is only very slightly different from the index  
$L(F,a,b)$ considered by Bourgain \rf{8}. $L(F;a,b) = \inf \{\eta <\omega_1  
\mid$ there exists a transfinite increasing sequence of open sets  
$(G_\alpha)_{\alpha \le\eta}$ with $G_0 =\emptyset$, $G_\eta = K$,   
$G_{\alpha+1}\setminus G_\alpha$ is disjoint from either   
$\lbrack F\le a\rbrack$ or  
$\lbrack F\ge b\rbrack$ for all   
$\alpha <\eta$ and $G_\gamma =\bigcup_{\alpha <\gamma}  
G_\alpha$ if $\gamma\le \eta$ is a limit ordinal$\}$.  
In fact one can show that   
if $\alpha (F;a,b) = \eta +n$ where $\eta$ is a limit ordinal and  
$n\in\IN$, then $L(F,a,b)\in \{ \eta + 2n,\, \eta +2n-1\}$. In  
Proposition~2.3 we shall show that $\alpha (F;a,b)<\omega$ for all  
$a<b$ iff $\beta (F,\delta) <\omega$ for all $\delta >0$. We note that  
a more general result has subsequently been obtained in \rf{25}.  
Indeed if we define $\beta (F) = \sup \{\beta (F;\delta)\mid \delta>0\}$  
and $\alpha (F) = \sup \{ \alpha (F;a,b)\mid a<b$ rational$\}$ then  
Kechris and Louveau have shown that  
$\beta (F) \le \omega^\xi$ iff $\alpha (F) \le \omega^\xi$.   

Also we note that the following result follows from \rf{8}. Let $X$  
be a separable Banach space not containing $\ell_1$. Let $K= Ba(X^*)$  
in its weak* topology. Then   
$$\sup \bigl\{ \beta (x^{**}|_K) : x^{**} \in X^{**}\bigr\} <\omega_1\ .$$  

\proclaim Proposition 2.3.  
Let $F: K\to \IR$ be a bounded function on the compact metric space $K$.  
The following are equivalent  
\smallskip  
\iitem{1)} $F\in B_{1/2} (K)$.  
\iitem{2)} $F$ is the uniform limit of $\D$-simple functions on $K$.  
\iitem{3)} For $a<b$, $\lbrack F\le a\rbrack$   
and $\lbrack F\ge b\rbrack $ may be separated by  
disjoint sets in $\D$.  
\iitem{4)} $\beta (F)\le \omega$.  
\iitem{5)} $\alpha (F; a,b) <\omega$ for all $a<b$.  
\smallskip\par  

\proof   

$4) \Rightarrow 5)$. This follows from the elementary observation that  
for all ordinals $\alpha$ and reals $a<b$, $K_\alpha (F;a,b) \subseteq  
K_\alpha (F,b-a)$, and the fact that 4) holds if and only if   
$\beta (F,\delta) <\omega$ for all $\delta >0$.  

$5)\Rightarrow 3)$.   
Let $K_i = K_i (F;a,b)$. Thus $K=K_0 \supseteq K_1\supseteq \cdots\supseteq  
K_n = \emptyset$ where $n=\alpha (F;a,b)$. Let   
$$D= \bigcup_{i=1}^n   
\ov{(F\le a\cap K_{i-1})}\setminus \ov{([F\ge b]\cap K_{i-1})}\in \D\ .$$  
Since $K_i= \ov{(\lbrack F\le a\rbrack\cap K_{i-1})}\cap   
\ov{(\lbrack F\ge b\rbrack \cap K_{i-1})}$,  
$$\eqalign{ D & = \bigcup_{i=1}^n \bigl( \, \ov{\lbrack F\le a\rbrack   
\cap K_{i-1}}  
	\setminus K_i\bigr)\cr  
\noalign{\vskip6pt}  
&\supseteq \bigcup_{i=1}^n \Big\lbrack \bigl(\lbrack F\le a\rbrack    
\cap K_{i-1}\bigr)  
	\setminus K_i\Big\rbrack\cr  
\noalign{\vskip6pt}  
&= \bigcup_{i=1}^n \bigl( \lbrack F\le a\rbrack   
\cap (K_{i-1}\setminus K_i)\bigr)  
	= \lbrack F\le a\rbrack \ .\cr}$$  
Since $K_{i-1}$ is closed,  
$$\eqalign{ D & \subseteq \bigcup_{i=1}^n \bigl( K_{i-1} \setminus  
	\ov{\lbrack F\ge b\rbrack \cap K_{i-1}}\ \bigr)\cr  
\noalign{\vskip6pt}  
& \subseteq \bigcup_{i=1}^n \Big\lbrack K_{i-1}   
\setminus \bigl( \lbrack F\ge b\rbrack   
	\cap K_{i-1}\bigr) \Big\rbrack\cr  
\noalign{\vskip6pt}  
& =\bigcup_{i=1}^n \bigl( K_{i-1}\setminus \lbrack F\ge b\rbrack \bigr)  
	= K\setminus \lbrack F\ge b\rbrack\ .\cr}$$	   

$3)\Rightarrow 2)$. This is a standard exercise in real analysis.  

$2)\Rightarrow 1)$. Since every $\D$-simple function can be expressed  
in the form $\sum_{i=1}^k a_i \bone_{L_i}$ where the $L_i$'s are closed sets  
and $DBSC(K)$ is a linear space it suffices to   
recall that $\bone_L\in DBSC(K)$  
whenever $L$ is closed. In fact $\bone_L$ is upper semicontinuous.  

$1)\Rightarrow 4)$. Let $F$ be the uniform limit of $(F_n)\subseteq DBSC(K)$. 
For $\delta >0$ and $n$ sufficiently large, $\beta (F,2\delta)  
\le \beta (F_n,\delta)$ and thus is suffices to prove that for  
$G\in DBSC(K)$, $\beta (G,\delta) <\omega$ for $\delta >0$. This is 
immediate from the following 

\proclaim Lemma 2.4. 
If $m\in\IN$, $\delta >0$ and $G: K\to \IR$ is such that $K_m(G,\delta) 
\ne \emptyset$, then $|G|_D \ge m\delta /4$.\par 

\proof Let $(g_n)\subseteq C(K)$ converge pointwise to $G$. It suffices 
to show that there exist integers $n_1<n_2 <\cdots < n_{m+1}$ and 
$k\in K$ such that $|g_{n_{i+1}} (k) - g_{n_i} (k)| > \delta /4$ for   
$1\le i\le m$.   

Let $n_1 =1$, $k_0 \in K_m (G,\delta)$ and let $U_0$ be a neighborhood   
of $k_0$ for which $\osc_{U_0} g_{n_1} < \delta /8$. Choose $k_0^1$   
and $k_0^2$ in $U_0\,\cap\, K_{m-1} (G,\delta)$ with $G(k_0^1) - G(k_0^2)   
> 3\delta/4$. Then choose $n_2>n_1$ such that $g_{n_2} (k_0^1) - g_{n_2}   
(k_0^2) > 3\delta/4$. Thus there is a nonempty neighborhood $U_1 \subset   
U_0$ of either $k_0^1$ or $k_0^2$ such that for $k\in U_1$, $|g_{n_2}   
(k) - g_{n_1}(k)| >\delta/4$.    

Similarly we can find a neighborhood $U_2 \subseteq U_1$ of a point in   
$K_{m-1} (G,\delta)$ and $n_3 >n_2$ so that for $k\in U_2$, $|g_{n_3}   
(k) - g_{n_2}(k)| > \delta/4$, etc.~\qed   

\demo Remarks 2.5. 1. Of course by using a bit more care one can show that   
$|G|_D \ge m\delta/2$ whenever $K_m (G,\delta) \ne \emptyset$.   

2. Following \rf{25} we say that for $F\in B_1(K)$, $F\in B_1^\xi(K)$ iff   
$\beta (F) \le \omega^\xi$. Thus $B_{1/2} (K) \equiv B_1^1 (K)$ by   
Proposition~2.3, a result also observed in \rf{25}.   

3. We do not yet have an index characterization of $B_{1/4}(K)$, however   
we have a necessary condition (which  may be sufficient). To describe   
this we first must generalize our index above. Let $F:K\to \IR$    
and let $(\delta_i)_{i=1}^\infty$ be  positive numbers.    
Set $K_0 (F,(\delta_i))   
= K$ and for $0\le i$ set   
$$K_{i+1} \bigl( F,(\delta_j)\bigr) = \bigl\{ k\in K_i \bigl( F,   
(\delta_j)\bigr) \mid \osc_{K_i(F,(\delta_j))}   
(F,k) \ge\delta_{i+1}\bigr\}\ .$$   

\proclaim Proposition 2.6.   
Let $F\in B_{1/4}(K)$. Then there exists an $M <\infty$ so that if $K_n   
(F,(\delta_i))\ne \emptyset$, then $\sum_{i=1}^n \delta_i \le M$.\par   

\proof   
Let $F$ be the uniform limit of $(G_n)$ with $|G_n|_D \le C <\infty$   
for all $n$. Suppose that $K_n (F,(\delta_i))\ne \emptyset$ for   
some sequence $(\delta_i)_{i=1}^\infty \subseteq \IR^+$. Since   
$K_n(F,(\delta_i))\subseteq K_n (G_m,(\delta_i/2))$ for large $m$,   
the latter set is non-empty as well. The proof of Lemma~2.4 yields   
$$\left\{   
\eqalign{&\hbox{If $G: K\to \IR$ and   
$(\delta_i)_{i=1}^\infty \subseteq \IR^+$   
	is such that $K_n(G,(\delta_i))\ne \emptyset$,}\cr   
&\hbox{then $|G|_D \ge 4^{-1} \sum\limits_{i=1}^n \delta_i$.} \cr}   
\right.\leqno(2.1)$$   

Thus by (2.1) we have, for large $m$, $C\ge |G_m|_D \ge 4^{-1} \sum_{i=1}^n   
\delta_i$    

\noindent and so $\sum_{i=1}^n \delta_i \le 4C$.~\qed   
\medskip   

We shall explore in greater detail in \S3 and \S8 some questions   
related to the problem of an index characterization of Baire-1/4. The   
following proposition gives a sufficient index criterion for a function   
to be Baire-1/4.  It also shows (via Proposition~2.3) that if   
$F\in B_{1/2} (K)\setminus B_{1/4}(K)$, then $\beta (F) = \omega$.   

\proclaim Proposition 2.7.   
Let $F\in B_1 (K)$. If $\beta (F)<\omega$, then $F\in B_{1/4} (K)$.\par   

\proof   
Without loss of generality let $F: K\to \lb 0,1\rb $ with $\beta (F) \le n$.   
Thus $\alpha (F;a,b)\le n$ for all $a<b$. It follows from the   
proof of $5)\Rightarrow 3)$ in Proposition~2.3 that for all $0<a<b<1$   
there exists a $D\in \D$ with $|\bone_D|_D \le 2n$,  $\lb F\le a\rb \subseteq   
D$ and $\lb F\ge b\rb \,\cap\, D=\emptyset$. Thus for all $m< \infty$ there   
exist sets $D_1\supseteq D_2\supseteq \cdots \supseteq D_m$ in $\D$   
with $\lb F\ge i/m\rb \subseteq D_i$, $\lb F\le (i-1)/m\rb \,\cap\, D_i   
=\emptyset$ and $|\bone_{D_i}|_D \le 2n$ for $i\le m$. In particular   
if $G=\sum_{i=1}^m  m^{-1} \bone_{D_i}$, then $\Vert F-G\Vert_\infty   
\le  m^{-1}$ and $|G|_D \le 2n$.~\qed   
\medskip   

The following proposition   
is related to work of A.~Sersouri \rf{39}.   
It is of interest to us because it shows that a separable   
Banach space $X$ can have functions of large index in $X^{**}$ and yet   
be quite nice.  In fact it shows there are Baire-1 functions of   
arbitrarily large index which strictly govern the class of    
quasireflexive (order~1) Banach spaces.   
Our proof was motivated by discussions with A.~Pe{\l}czy\'nski.   

\proclaim Proposition 2.8.   
For all $\gamma <\omega_1$ there exists a quasireflexive (of order $1$)   
Banach space $Q_\gamma$ such that $Q_\gamma^{**} = Q_\gamma \oplus   
\langle F_\gamma\rangle$ where $\beta (F_\gamma) > \gamma$.\par   

\noindent   
(The index $\beta (F_\gamma)$ is computed with respect to $Ba (Q_\gamma^*)$.)   

\demo Remark 2.9.   
In \S6 we shall show the existence of a quasireflexive space whose new   
functional (in the second dual) is Baire-1/4.   

\noindent {\it Proof of Proposition 2.8.}   
We use interpolation, namely the method of \rf{12}. (This has also been   
used in \rf{19} in a slightly different manner to produce a    
quasireflexive  space from a weak* convergent sequence.)   

To begin let $\gamma < \omega_1$ be any ordinal and choose a compact   
metric space $K$ containing an ambiguous set $A_\gamma$ with $\alpha    
(\bone_{A_\gamma} ; \frac14,\frac34) >\gamma$.   
(For example $\bone_{A_\alpha}$   
could be taken to be one of the functions $F_\delta$ described in   
\S5 with $\delta >\omega^\gamma+$.) Choose a sequence $(f_n)\subseteq   
Ba (C(K))$ converging pointwise to $\bone_{A_\gamma}$ such that   
$(\bone_{A_\gamma},f_1,f_2,\ldots)$    
is basic in $C(K)^{**}$. Let $W$ be the closed   
convex hull of $\{ \pm f_n\}_{n=1}^\infty $ in $C(K)$. Let $Q_\gamma$   
be the Banach space obtained from $W\subseteq Ba(C(K))$ by    
$\lb$DFJP$\rb$-interpolation.   
Thus for all $n\in\IN$, $\Vert\cdot\Vert_n$ is the gauge of $U_n =   
2^nW +2^{-n}Ba(C(K))$, and $Q_\gamma = \{ x\in C(K) :   
\trivert x\trivert \equiv (\sum_n \Vert x\Vert_n^2)^{1/2} <\infty\}$.   
Following the notation of \rf{12}, we let $C= Ba(Q_\gamma) =   
\{ x\in C(K) :\trivert x\trivert \le 1\}$ and let $j: Q_\gamma \to   
C(K)$ be the natural semiembedding.   

We first observe that $Q_\gamma$ is quasireflexive of order $1$. Indeed   
it is easy to check that $\widetilde W$, the weak* closure of $W$ in   
$C(K)^{**}$ is just   
$$\widetilde W = \biggl\{ \sum_{i=1}^\infty a_if_i+a_\infty \bone_{A_\gamma}   
: |a_\infty| +\sum_{i=1}^\infty |a_i|\le 1\biggr\}\ .$$  
Furthermore $\widetilde C \subseteq [\widetilde W]$ (\rf{12}, Lemma~1(v))  
which has the basis $(\bone_{A_\gamma} , f_1,f_2,\ldots)$. Now  
$j^{**} : Q_\gamma^{**} \to C(K)^{**}$ is one-to-one and $(j^{**})^{-1}  
(C(K)) = Q_\gamma$ (Lemma~1(iii)). Thus if $F_\gamma\in Q_\gamma^{**}$  
satisfies $j^{**} F_\gamma =\bone_{A_\gamma}$, then $Q_\gamma^{**} = Q_\gamma  
\oplus \langle F_\gamma\rangle$. Of course $F_\gamma$ must be the weak*  
limit of $(j^{-1}(f_n))_n$ in $Q_\gamma^{**}$.  

It remains to show that $\beta (F_\gamma)\ge \gamma$. We shall prove  
$$\ov{\alpha} (F_\gamma; \frac14,\frac34) \ge \alpha (\bone_{A_\gamma} ;  
\frac14,\frac34) \leqno(2.2)$$  
where $\ov{\beta}$ is the index computed with respect to $F_\gamma \in  
B_1 (3Ba(Q_\gamma^*))$. Since $\beta (F_\gamma) \ge \alpha  (F_\gamma;  
\frac1{12},\frac14) \ge \ov{\alpha} (F_\gamma; \frac14,\frac34)$,  
the result follows.  

Since $\Vert j\Vert \le3$, if $K_0 = 3Ba (Q_\gamma^*)$ and $H_0= Ba   
(C(K)^*)$, then $j^* H_0\subseteq K_0$. More generally if  
$K_{\beta +1} = \{ y^* \in K_\beta \mid$ for all non-empty relative  
weak* neighborhoods $U$ of $y^*$ in $K_\beta$ there exists $y_1^*,y_2^*  
\in U$ with $F_\gamma (y_1^*)\ge \frac34$ and $F_\gamma (y_2^*)\le\frac14\}$  
and $H_{\beta +1} $ is defined similarly in terms of $\bone_{A_\gamma}$, 
then $j^* H_{\beta +1} \subseteq K_{\beta +1}$ for all $\beta$,  
since $j^*$ is $\omega^*$-continuous and $F_\gamma (j^* x^*) = (j^{**}  
F_\gamma)x^* =\bone_{A_\gamma} (x^*)$. This proves~(2.2).~\qed 

\beginsection{3. Theorem B.} 

For the proof of Theorem (B) (a) we need a lemma. Recall that a  collection 
of pairs of subsets of $K$, $(A_i,B_i)_{i=1}^n$,   
is said to be ({\it Boolean\/}) 
{\it independent\/} if for all $I \subseteq \{ 1,\ldots,n\}$,  
$\bigcap_{i\in I} A_i\ \cap\ \bigcap_{i\notin I} B_i \ne \emptyset$.  

\proclaim Lemma 3.1.  
Let $F: K\to \IR$ be the pointwise limit of $(f_n) \subseteq C(K)$.  
If $K_m (F;a,b)\ne\emptyset$ for some $m\in \IN$ and $a<b$, then  
for $a<a'<b'<b$ there exists a subsequence  $(f'_n)$ of $(f_n)$ so that  
if $n_1<\cdots < n_m$, then $(A'_{n_i},B'_{n_i})_{i=1}^m$ are  
independent where $A'_{n_i}=\lb f'_{n_i}\le a'\rb$   
and $B'_{n_i} = \lb f'_{n_i}  \ge b'\rb $. \par

\proof  
The proof is similar to that of Lemma~2.4 and is actually a local version  
of the proof of the main result of \rf{35} (see \rf{8} for a more general  
discussion of the consequences of $K_\beta (F;a,b)\ne\emptyset$).  

We first show how to choose a finite subsequence $(f_{n_i})_{n=1}^m$ of  
$(f_n)$ so that $(A_{n_i},B_{n_i})_{i=1}^m$ is independent, where  
$A_{n_i} = \lb f_{n_i} \le a'\rb$ and $B_{n_i} = \lb f_{n_i}\ge b'\rb$. Let   
$k_\phi \in K_m (F;a,b)$. Thus there  exist $k_0$ and $k_1$ in  
$K_{m-1} (F;a,b)$ with $F(k_0) \le a$ and $F(k_1) \ge b$. Choose  
$n_1$ and neighborhoods $U_0$ and $U_1$ of $k_0$ and $k_1$, respectively,  
so that $f_{n_1} < a'$ on $U_0$ and $f_{n_1} >b'$ on $U_1$. Let  
$k_{\varep_1,\varep_2}\in U_{\varep_1}\,\cap\, K_{m-2} (F;a,b)$  
for $\varep_1,\varep_2\in \{ 0,1\}$  with $F(k_{\varep_1,0})\le a$  
and $F(k_{\varep_1,1})\ge b$ for $\varep_1 \in \{ 0,1\}$.  
Choose $n_2>n_1$ and neighborhoods $U_{\varep_1,\varep_2}    
\subseteq U_{\varep_1}$  
of $k_{\varep_1,\varep_2}$ so that $f_{n_2} <a'$ on $U_{\varep_1,0}$ and  
$f_{n_2} >b'$ on $U_{\varep_1,1}$ (for $\varep_1,\varep_2\in\{ 0,1\}$).  
Continue up to $f_{n_m} $. The sets $(A_{n_i},B_{n_i})_1^m$ are  
then independent since for $I\subseteq \{ 1,\ldots,m\}$, $\bigcap_{i\in I}  
A_{n_i}\cap \bigcap_{i\notin I} B_{n_i} \supseteq U_{\varep_1,\ldots,  
\varep_m} \ne \emptyset$ where $\varep_i =0$ if $i\in I$ and   
$\varep_i = 1$ if $i\notin I$.  

Now the existence of an infinite subsequence $(f'_n)$ satisfying the  
conclusion of 3.1 follows immediately from Ramsey's theorem. Indeed,  
by the latter, there exists $(f'_n)$ a subsequence of $(f_n)$ so that  
$(f'_n)$ satisfies the conclusion, or such that for {\it all\/}  
$n_1<\cdots < n_m$, $(A'_{n_i}, B'_{n_i})_{i=1}^m$ {\it is  not\/}  
independent. But we have proved that the second alternative is impossible.  
\medskip  

\noindent {\it Proof of Theorem B(a)}.  
$(f_n)$ is a bounded sequence in $C(K)$ converging pointwise to $F\notin  
B_{1/2} (K)$. By Proposition~2.3 there exists $a<b$ so that $K_m (F;a,b)\ne   
\emptyset$ for all $m\in \IN$. By passing to a subsequence we may assume  
$(f_n)$ has a spreading model.  
Furthermore by Lemma~3.1, passing to subsequences and diagonalization  
we may assume that for some  
$a<a'<b'<b$, $(A_{n_i},B_{n_i})_{i=1}^m$ is independent  
whenever $m\le n_1<n_2< \cdots < n_m$ and $A_{n_i} = \lb f_{n_i} \le a'\rb$,  
$B_{n_i} = \lb f_{n_i} \ge b'\rb$.   
By Proposition~4 of \rf{36} it follows that  
there exists $C<\infty$ so that $(f_{n_i})_{i=1}^m$ is $C$-equivalent  
to the unit vector basis of $\ell_1^m$   
whenever $m\le n_1<\cdots < n_m$.~\qed  
\medskip  

The proof of Theorem B(b) will require a more precise version of Theorem~A(b)  
and the following elementary lemma (which follows easily from the  
Hahn-Banach theorem). If $C$ is a subset of a Banach space $X$,   
$\widetilde C$ denotes the  $w$*-closure of $C$ in $X^{**}$.   

\proclaim Lemma 3.3.  
Let $C$ and $D$ be convex subsets of $X$. Then $md (C,D) = md (\widetilde C,  
\widetilde D)$. By $md(C,D)$ we mean the minimum distance,  
$$\inf \bigl\{ \Vert c-d\Vert \mid c\in C\ ,\ d\in D\bigr\}\ .$$\par  

The variant of Theorem A(b) which we need is   

\proclaim Lemma 3.4.  
Let $F: K\to \IR$ be bounded and let $(f_n) \subseteq C(K)$  
converge pointwise to $F$ with $\sum_{n=0}^\infty |f_{n+1} (k) - f_n (k)|  
\le M$ for all $k\in K$ $(f_0\equiv 0)$. Suppose $\osc (F,k) >\delta$  
for some $\delta >0$. Then there exists a subsequence $(f'_n)$ of $(f_n)$  
which is $C= C(M,\delta)$ equivalent to the summing basis.\par  

Let $F\in B_1(K)\setminus C(K)$. It is evident that if $F$ strictly  
governs $\{ c_0\}$, then $F\in DBSC(K)$. The next result shows that  
the converse is true.  

\proclaim Corollary 3.5.  
Let $F\in DBSC(K)$ and let $(f_n)$, $M$ and $\delta$ be as in the hypothesis   
of Lemma~3.4. Let $(g_n)\subseteq C(K)$ converge pointwise to $F$  
with $\sup_n \Vert g_n\Vert_\infty < \infty$. Then there exists $(h_n)$,  
a convex block subsequence of $(g_n)$, which is $C(M,\delta)$-equivalent  
to the summing basis.\par  

The proof is straightforward from Lemmas 3.3 and 3.4.  
\medskip  

\noindent {\it Proof of Theorem B(b)}. Let $F\in B_{1/4} (K)\setminus C(K)$  
and let $(f_n) \subseteq C(K)$ be a bounded sequence converging pointwise  
to $F$. Choose $(F_n) \subseteq DBSC(K)$ which converges uniformly to $F$  
so that $\sup_n|F_n|_D < M<\infty$. For each $n\in \IN$, choose  
$(f_i^n)_{i=0}^\infty \subseteq C(K)$, $f_0^n\equiv 0$, which  
converges pointwise to $F_n$ and satisfies $\sum_{i=0}^\infty |f_{i+1}^n  
(k) - f_i^n(k)|\le M$ for $k\in K$.  

Since $F\notin C(K)$ we may assume there exists $\delta >0$ so that for all  
$n$, $\osc_K (F_n,k) >\delta >0$ for some $k\in K$. Thus, by Lemma~3.4,  
we may suppose for all $n$, $(f_i^n)_{i=1}^\infty$ is $C=   
C(M,\delta)$-equivalent to the summing  basis. We may also assume   
$\Vert F_n-F_{n+1}\Vert_\infty <\varep_n$  where $\varep_n\downarrow 0$  
and for all $n\in \IN$, $\sum_{i=n+1}^\infty \varep_i <\varep_n$.  

By induction and Lemma 3.3 we may replace each sequence $(f_i^n)_{i=1}^\infty$  
by a convex block subsequence $(g_i^n)_{i=1}^\infty$ such that for $n>1$,  
$$\left\{  
\eqalign{&\hbox{there exists a convex block subsequence   
	$(h_i^n)_{i=1}^\infty$ of $(g_i^{n-1})_{i=1}^\infty$}\cr   
&\hbox{with $\Vert g_i^n-h_i^n\Vert_\infty <\varep_{n-1}$ for $i\in\IN$.}\cr}  
\right. \leqno(*)$$   

Let $(g_n^n)_{n=1}^\infty$ be the diagonal sequence. Clearly $(g_n^n)$   
converges pointwise to $F$. Also by $(*)$ for $n>k$,    
$md (g_n^n,co(g_j^k)_{j=1}^\infty) < \sum_{j=k}^n \varep_j <   
\varep_{k -1}$. In fact for $k$ fixed, there exists a convex block   
subsequence $(d_n^k)_{n>k}$ of $(g_j^k)_{j=1}^\infty$ with    
$\Vert g_n^n-d_n^k\Vert_\infty <\varep_{k-1}$ for $n>k$.   
Thus for any $k$, $(g_n^n)_{n>k}$ is an $\varep_{k-1}$-perturbation   
of a sequence $(d_n^k)_{n>k}$ which is $C'$-equivalent to the summing   
basis where $C'$ depends solely on $C$.   

By Lemma 3.3 applied to $co (f_n)$ and $co(g_n^n)$, there are convex   
block subsequences $(g_n)$ of $(f_n)$ and $(\bar g_n)$ of $(g_n^n)$   
with $\Vert g_n-\bar g_n\Vert_\infty \to 0$. Since $(\bar g_n)_{n>i}$   
is an $\varep_{i-1}$-perturbation of a sequence which is   
$C'$-equivalent to the summing basis, $(\bar g_n)$ and hence $(g_n)$   
has a subsequence which has spreading model equivalent to the summing   
basis.~\qed   

\demo Remark 3.6.   
The constant of equivalence of the spreading model of $(g_n)$ with the   
summing basis depends solely upon $\sup_{k\in K} \osc_K(F,k)$ and   
$|F|_{1/4}$.\par   

Our next theorem is a converse to Theorem B(a).   

\proclaim Theorem 3.7.   
Let $F\in B_1 (K)$. Assume that whenever $(f_n) \subseteq C(K)$ is a    
uniformly bounded sequence converging pointwise to $F$, then any spreading   
model of $(f_n)$ is equivalent to the unit vector basis of $\ell_1$.   
Then $F\notin B_{1/2} (K)$.\par   

\proclaim Lemma 3.8.   
Let $F\in B_{1/2}(K)\setminus C(K)$, $\Vert F\Vert_\infty \le1$.   
Then there exists $(f_n) \subseteq C(K)$ converging pointwise to $F$   
with spreading model $(e_n)$ and a function $M : \IR^+ \to \IR^+$   
satisfying   
$$\Vert \sum a_ne_n\Vert \le M(\varep) \Vert \sum a_ns_n\Vert   
+ \varep \sum |a_n|\leqno(3.1)$$   
for all $(a_n) \subseteq \IR$ and $\varep >0$.\par   

\proof Let  $(g_n) \subseteq Ba (C(K))$ converge pointwise to $F$ and let   
$\varep_n\downarrow 0$. By the proof of Theorem~B(b) we can choose $(f_n)$,   
a convex block subsequence of $(g_n)$ such that for all $m$,    
$(f_n)_{n=m}^\infty$   
is an $\varep_m$-perturbation of a sequence which is   
$M(\varep_m,F)$-equivalent   
to the summing basis.~\qed  
\medskip  

\noindent {\it Proof of Theorem 3.7.}  
This is immediate from Lemma 3.8, since if $(e_n)$ satisfies (3.1), then   
$$\lim_{n\to\infty}\ {1\over n} \Big\Vert \sum_{i=1}^n (-1)^i e_i  
\Big\Vert =0\ .$$  
In particular $(e_i)$ is not equivalent to the unit vector basis of  
$\ell_1$.~\qed  
\medskip  

The proof of Theorem 3.7 combined with Theorem B(a) yields the   
following result.  {\sl Let $F\in B_1(K)$. Then $F\notin B_{1/2}(K)$  
if and only if there exists $(f_n) \subseteq C(K)$, a uniformly  
bounded sequence converging pointwise to $F$, so that if $(g_n)$ is a  
convex block subsequence of $(f_n)$, then some subsequence of $(g_n)$  
has the unit vector basis of $\ell^1$ as a spreading model\/}.  

We do not know if the converse to Theorem B(b) is valid.  

\demo Problem 3.9.  
Let $F\in B_1 (K)$ and $C <\infty$ be such that whenever $(f_n)$ is a         
uniformly bounded sequence in $C(K)$ converging pointwise to $F$, then there   
exists $(g_n)$, a convex block subsequence of $(f_n)$ with spreading  
model $C$-equivalent to the summing basis. Is $F\in B_{1/4}(K)$?\par  

We now turn to the Banach space implications of Theorem B. Let $K$ be  
compact metric and let $X$ be a closed subspace of $C(K)$. For example,  
$K$ could be $Ba (X^*)$ but we do not require this. $X^{**}$ is   
naturally isometric to $X^{\bot\bot}\subseteq C(K)^{**}$. In this setting  
it can be shown (see \rf{35}) that if $B_1(X) = \{ x^{**} \in X^{**}:$  
there exists $(x_n)\subseteq X$ with $(x_n)$ converging weak* in $X^{**}$  
to $x^{**}\}$, then $B_1(X) \subseteq B_1(C(K))$ and $B_1 (C(K))$ is  
naturally identified with $B_1(K)$.  

\proclaim Corollary 3.10.  
Let $K$ be compact metric and let $X$ be a closed subspace of $C(K)$.  
\smallskip  
\iitem{a)} If $X^{**}\,\cap\,\lb B_1(K)\setminus B_{1/2}(K)\rb\ne\emptyset$,  
then $X$ contains a basic sequence with spreading model equivalent to the  
unit vector basis of $\ell_1$.  
\iitem{b)} If $\lb X^{**}\,\cap\,B_{1/4}(K)\rb \setminus X\ne \emptyset$  
then $X$ contains a basic sequence with spreading model equivalent to  
the summing basis.  
\smallskip\par  

\demo Remark 3.11.  
This corollary has immediate purely local consequences. Thus if $X$  
and $K$ are as above and $X$ does not contain $\ell_n^\infty$'s uniformly,  
then $X^{**} \cap B_{1/4} (K)\subset X$. Moreover if $X$ is $B$-convex,  
{\it i.e.}, does not contain $\ell^1_n$'s uniformly, then  
$X^{**}\setminus X\subset B_{1/2} (K) \setminus B_{1/4}(K)$.\par

\beginsection{4. $DSC(K)$.}  

\proclaim Theorem 4.1.  
Let $K$ be compact metric and let $F\in DSC(K)\setminus C(K)$. Then  
$F$ governs $\{ c_0\}$.\par  

\demo Remark 4.2.  
If $X$ is a separable Banach space, $K= Ba  (X^*)$ in its weak* topology  
and $F\in X^{**}$, then if $F\in DSC(K)$, $F\in DBSC(K)$ (and hence  
for such functions Theorem~4.1 follows from Theorem~A). To see this  
assertion, first choose $(f_n)$ uniformly bounded in $C(K)$ so that  
$f_n\to F$ pointwise and $\sum_{n=1}^\infty |f_{n+1}(k) -f_n(k)| <\infty$  
for all $k\in K$. Now since $F\in B_1(X)$, we may choose $(g_j)$ a   
convex block subsequence of $(f_j)$ and $(x_j)$ a sequence in $X$ with  
$\| g_j-x_j\| <  2^{-j} $ for all $j$. But then it follows that   
$x_j \to F$ pointwise and moreover $\sum_{j=1}^\infty |x_{j+1} (k) -x_j(k)|  
< \infty$ for all $k\in K$. Thus by the uniform boundedness principle,  
$$\sup_{k\in K} \sum_{j=1}^\infty | x_{j+1} (k) - x_j(k) | <\infty\ ,$$  
so $F\in DBSC(K)$.  

Theorem 4.1 follows from the stronger result of Elton \rf{13}  
which was motivated by work of Fonf \rf{16}.  

\proclaim Theorem.  \rf{13}.    
Let $X$ be a Banach space and let $\E$ be the set of extreme points  
of $Ba(X^*)$. Let $(x_i)$ be a normalized basic sequence in $X$ such that  
$\sum_{i=1}^\infty |x^*(x_i)| <\infty$ for all $x^*\in \E$. Then   
$c_0 \hookrightarrow \lb (x_i)\rb$.\par  

Theorem 4.1. can be phrased in  this way provided $\E$ is replaced  
by $\ov{\E}$. However we wish to present a separate proof of our weaker  
result which seems to be of interest in its own right. The main step is  
given by the following lemma. If $Y$ is a subspace of $C(K)$, $U\subseteq K$  
and $r>0$, we say {\it $U$ $r$-norms $Y$\/} if $\Vert y\big|_U\Vert_\infty\ge   
r\Vert y\Vert$ for all  $y\in  Y$.  

\proclaim Lemma 4.3.  
Let $L$ be a compact metric space and let $(f_i)$ be a normalized basic  
sequence in $C(L)$. If $c_0\not\hookrightarrow [(f_i)]$, then there exists  
a nonempty compact set $K\subseteq L$ and a normalized block basis $(g_i)$  
of $(f_i)$ so that  
$$\left\{ \eqalign{&\hbox{for any nonempty relatively open subset $U$ of $K$   
	there are an}\cr  
&\hbox{$r>0$ and an $n_0\in \IN$ such that $U$ $r$-norms   
	$[(g_n)_{n=n_0}^\infty]$.}\cr}\right.\leqno(4.1)$$\par  

\demo Remark 4.4.  
It can be deduced from \rf{36}    
that $[(x_n)]$ contains an isomorph of $\ell_1$ iff  
there exists a compact  set $K\subseteq L$ such that (4.1) holds for some  
fixed $r>0$ independent of $U$.\par  

\noindent {\it Proof of Lemma 4.3.}  
Let $(U_m)_{m=1}^\infty$ be a base of open sets for $L$. We inductively  
construct for each $m$ a normalized  block basis $(f_i^m)_{i=1}^\infty$  
of $(f_i)$ and a certain subsequence $M$ of $\IN$.  

Let $(f_i^0) = (f_i)$ and suppose $(f_i^m)_{i=1}^\infty$ has been chosen.  
There are two possibilities.   
\smallskip  
\iitem{(i)} There is a normalized block basis $(g_i)$ of   
$(f_i^m)_{i=1}^\infty$  
with $\Vert g_i\big|_{U_m}\Vert_\infty \to 0$ as $i\to \infty$.  
\iitem{(ii)} There exists no such sequence.  
\smallskip  
\noindent If (i) holds, choose $(f_i^{m+1})_{i=1}^\infty$ to be a normalized  
block basis of $(f_i^m)_{i=1}^\infty $ with  
$$\Vert f\big|_{U_m}\Vert_\infty < 2^{-k} \Vert f\Vert_\infty  
\ \hbox{ for all }\ f\in \big\lb (f_i^{m+1})_{i=k}^\infty\big\rb  
\leqno(4.2)$$  
and put $m$ in $M$. If (ii) holds let   
$(f_i^{m+1})_{i=1}^\infty = (f_i^m)_{i=1}  
^\infty$ and put $m$ in $\IN \setminus M$.   
Let $K= L\setminus \bigcup_{m\in M}  
U_m$ and for all $n\in M$ let $g_n = f_{n+1}^{n+1}$. We may assume $M$ is  
infinite or else the conclusion of the lemma is satisfied with $K=L$  
and $g_i= f_i^m$ ($m=\max M$ or 0 if $M=\emptyset$).  

First we check that $K\ne\emptyset$. If $K= \emptyset$, then $L\subseteq  
\bigcup_{n\in M} U_n$. By compactness there exists $n_1 \in M$ so that  
$L\subseteq \bigcup_{n\in M\ ,\ n\le n_1} U_n$.  But then since  
$\Vert g_{n_1} \big|_{U_n}\Vert_\infty <2^{-(n_1+1)}$ for $n\in M$  
with $n\le n_1$, we have $\Vert g_{n_1}\Vert_\infty < 1$, a contradiction.  

We claim that $K$ and $(g_n)$ satisfy (4.1). If not there exist  
$(h_n)$, a normalized block basis of $(g_n)$ and a $U_m$ such that  
$K\cap U_m \ne \emptyset$ and so   
$m\notin M$ yet $\Vert  h_i\big|_{K\cap \ov{U}_m}  
\Vert < 2^{-i}$ for all $i$. Indeed there must exist $m'\in M$  
with $K\cap U_{m'} \ne \emptyset$ and $(h_i)$, a normalized block  
basis of $(g_n)$, with $\Vert h_i\big|_{K\cap U_{m'}}\Vert < 2^{-i}$.  
Then choose $m\in\IN$ so that $\ov{U}_m \subseteq U_{m'}$ and $K\cap U_m  
\ne\emptyset$. Let $j_0=m$ and if $j_i$ is defined choose $j_{i+1} > j_i$  
so that  
$$\ov{U}_m \cap\lb h_{j_i} \ge 2^{-i}\rb \subseteq    
\bigcup_{\scriptstyle n\in M  
\atop \scriptstyle n\le j_{i+1}} U_n\ .$$  
This can be done since $\ov{U}_m \cap \lb h_{j_i} \ge 2^{-i}\rb \subseteq  
\ov{U}_m \cap \lb h_{j_i} \ge 2^{-j_i}\rb \subseteq L\setminus  K=  
\bigcup_{n\in M} U_n$.  This completes the definition of $j_1,j_2,\ldots\ $.  
Now for $t\in U_m$, $|h_{j_i} (t)|\ge 2^{-i}$  
for at most one $i$. Indeed let $i_0$ be the first integer such that  
$|h_{j_{i_0}} (t)|\ge 2^{-i_0}$ (if such an $i_0$ exists). Then $t\in   
\bigcup_{n\in M\ ,\ n\le j_{i_0+1}} U_n$ and for $i>i_0$,  
$h_{j_i}$ is a normalized element in $\lb (g_j)_{j\ge j_i,\, j\in M}\rb     
= \lb (f_{j+1}^{j+1})_{j\ge j_i,\, j\in M}\rb  \subseteq \lb (f_p^{j_i+1} )    
_{p\ge j_i+1}\rb$. Thus if $t\in U_n$ with $n\le j_{i_0+1},n\in M$,  
then $h_{j_i} \in \lb (f_p^{n+1})_{p\ge j_i+1}\rb$ and so by (4.2),  
$|h_{j_i} (t)| \le \Vert h_{j_i}\big|_{U_n}\Vert < 2^{-j_i} \le 2^{-i}$.   

Thus $\sum_{i=1}^\infty |h_{j_i}(t)| \le 2$ for all $t\in \ov{U}_m$. Since  
$\ov{U}_m $ norms $\lb h_{j_i}\rb$, it follows from \rf7 that $c_0  
\hookrightarrow \lb h_{j_i}\rb$, a contradiction.~\qed  
\medskip  

\noindent {\it Proof of Theorem 4.1}.  
Let $(f_n)$ be a bounded sequence in $C(K)$ converging pointwise to $F$.  
By Lemma~3.3 and passing to a convex block subsequence of $(f_n)$,  
if necessary, we may suppose that $\sum_{n=1}^\infty |f_{n+1}(k)-f_n(k)|   
<\infty$  
for all $k\in K$.   
Also since $F\notin C(K)$, by passing to a subsequence $(f'_n)  
\subseteq (f_n)$ we may assume that $(h_n) \equiv (f'_{2n}-f'_{2n+1})$  
is a seminormalized basic sequence  satisfying $\sum_{n=1}^\infty  
|h_n(k)| < \infty$ for all $k\in K$.   
If $c_0\not\hookrightarrow  
\lb (h_n)\rb $, then by Lemma~4.3 there exist $(g_n)$,   
a normalized block basis  
of $(h_n)$, and a closed nonempty set $K_0\subseteq K$ satisfying (4.1)  
(with $K$ replaced by $K_0$).  

For $m\in\IN$ set $K_m = \{ k\in K_0 : \sum_{n=1}^\infty |g_n(k)| \le m\}$.  
Since $(g_n)$ is a normalized block basis of $(h_n)$, $\sum_{n=1}^\infty  
|g_n(k)| <\infty$ for all $k\in K$ and thus   
$\bigcup_{m=1}^\infty K_m = K_0$.  By the Baire category theorem  
there exists $m_0$ so that $K_{m_0}$ has nonempty interior $U$   
(relative to $K_0$). Choose $n_0$ and $r>0$ so that $U$ $r$-norms  
$\lb (g_n)_{n\ge n_0}\rb $. Since $\sum |g_n| \le m_0$ on $U$,  
$(g_n)$ is equivalent to the unit vector basis of $c_0$ \rf{7},   
a contradiction.~\qed  
\medskip  

A natural problem is to classify those functions $F\in B_1 (K)$ which  
govern $\{c_0\}$. We do not know how to  do this, but it is easy to see  
that  this class is strictly larger than $DSC(K)$.  

\demo Example 4.5.   
Let $L$ be a countable compact metric space, large enough so that there  
exists an $F\in B_1 (L) \setminus DBSC(L)$ (see Proposition~5.3). Choose a  
bounded sequence $(f_n) \subseteq C(L)$ which converges pointwise  
to $F$ and let $X = \lb (f_n)\rb$. $C(L)$ is $c_0$-saturated (every   
infinite dimensional subspace of $C(L)$ contains $c_0$ isomorphically)  
and thus   
$X$ is $c_0$-saturated. Thus $F$ governs $\{c_0\}$ by Lemma~3.3.  
Let $K = Ba (X^*)$. $F\notin DSC(K)$ or otherwise (Remark~4.2)   
$F\in DBSC(K)$ and hence $F\in DBSC(L)$.   
Using this example, it can be shown that if $K$ is any uncountable  
compact metric space, there exists an $F\in B_1(K)\setminus DSC (K)$  
which governs $\{ c_0\}$.  

\demo Question 4.6.  
Let $F\in B_1 (K)$. If $F$ governs $\{ c_0\}$ does there exist a bounded  
sequence $(f_n)\subseteq C(K)$ converging pointwise to $F$ and a $w^*$-closed  
set $L\subseteq Ba \lb (f_n)\rb^*$ such  that $L$ norms $\lb (f_n)\rb$ and  
$F\big|_L \in DSC(L)$? (Could $L$ be taken to be countable?)\par  

\demo Question 4.7.  
Let $F\in B_1 (K)$. Suppose there exists $(f_n) \subseteq C(K)$,  
a bounded sequence converging pointwise to $F$ and satisfying  
$\sum_{n=1}^\infty |f_{n+1} (k) -f_n (k) | <\infty$ for all $k$ in  
some residual set (complement of a first category set). Does $F$  
govern $\{ c_0\}$?\par  

We should also mention the following result of Bourgain which gives some  
global information about the class $DSC(K)$.  

\proclaim Proposition 4.8. \rf{10}\enspace  
Let $F\in DSC(K)\setminus C(K)$   
and let $(f_n)$ be a bounded sequence in $C(K)$ converging  
pointwise to $F$ with $\sum |f_{n+1}(k) - f(k)|<\infty$ for all $k\in K$.  
Then there exists a subsequence $(f'_n)$ of $(f_n)$ with $\lb (f'_n)\rb^*$   
separable.\par  

It follows that if $F\in DSC(K)\setminus C(K)$, then $F$ strictly governs  
the class $\C$ of infinite dimensional Banach spaces with separable duals.  
However we don't know that if $F$ governs $\{ c_0\}$, then  $F$ strictly  
governs $\C$. (A negative answer, of course, would give a negative answer  
to 4.6.)  

We give a somewhat different proof than that of \rf{10}.  

\proof We may assume that $\Vert f_n\Vert=1$ for all $n$.  
As mentioned in the introduction there exists a subsequence  
$(f'_n)$ of $(f_n)$ which is basic and $C_1$-dominates the summing basis for  
some $C_1<  \infty$. It follows that $(h_n)_1^\infty$ is seminormalized  
basic where $h_1 = f'_1$ and $h_n = f'_n - f'_{n-1}$ for $n>1$.   
$\lb$Indeed let $(a_i)_1^m$ be given and let $1\le n<m$  with   
$\Vert \sum_1^n a_ih_i\Vert =1$. $\sum_{i=1}^n a_i h_i = (a_1-a_2) f'_1+  
\cdots + (a_{n-1}-a_n) f'_{n-1} + a_nf'_n \equiv f + a_nf'_n$.  
If $\Vert f\Vert \ge \frac12$, then $\Vert \sum_1^m a_ih_i\Vert \ge C_2^{-1}  
\Vert f\Vert \ge 2^{-1} C_2^{-1}$ where $C_2$ is the basis constant  
of $(f'_n)$. Otherwise $|a_n| \ge \frac12$ and so   
$$\eqalign{ \Big\Vert \sum_{i=1}^m a_ih_i\Big\Vert  
&= \Big\Vert \sum_{i=1}^{m-1} (a_i -a_{i+1}) f'_i + a_mf'_m\Big\Vert\cr  
\noalign{\vskip6pt}  
&\ge (C_2+1)^{-1} \Big\Vert \sum_{i=n}^{m-1} (a_i - a_{i+1}) f'_i  
	+ a_m f'_m\Big\Vert\cr  
\noalign{\vskip6pt}  
&\ge (C_2 +1)^{-1} C_1^{-1} \Big\Vert \sum_{i=n}^{m-1} (a_i -a_{i+1})  
	s_i + a_m s_m\Big\Vert\cr  
\noalign{\vskip6pt}  
&\ge (C_2 +1)^{-1} C_1^{-1}\Big| \sum_{i=n}^{m-1} (a_i -a_{i+1}) +a_m\Big|\cr  
\noalign{\vskip6pt}  
&=(C_2+1)^{-1} C_1^{-1} |a_n| \ge 2^{-1}(C_2+1)^{-1}C_1^{-1}\ .\big\rb\cr}$$  

Also for $k\in K$, $\sum_{n=1}^\infty |h_n(k)| <\infty$.  
Thus $(h_n)$ is shrinking. Indeed if $(h_n)$ has basis constant $C$ and  
$g_n = \sum_{i=p_n+1}^{p_{n+1}} a_ih_i$ is a normalized block basis, then  
for $k\in K$  
$$\eqalign{ |g_n(k)| & \le \biggl(\max_{p_n+1\le i\le p_{n+1}}  
	|a_i| \biggr) \sum_{i=p_n+1}^{p_{n+1}} |h_i(k)| \cr  
\noalign{\vskip6pt}  
&\le (C+1) \min_i \Vert h_i\Vert^{-1} \sum_{i= p_n+1}^{p_{n+1}} |h_i(k)|\cr}$$  
which goes to 0 as $n\to\infty$.~\qed  
\medskip  

The following proposition characterizes the subclass $PS(K)$ of $DSC(K)$  
which was defined in \S1.  

\proclaim Proposition 4.9.  
Let $F\in B_1 (K)$. The following are equivalent.  
\smallskip  
\iitem{a)} $F\in PS (K)$.  
\iitem{b)} For all closed $L\subseteq K$, $F\big|_L$ is continuous  
on a relatively open dense subset of $L$.  
\iitem{c)} There exists $\eta <w_1$ and a family $(K_\alpha)_{\alpha\le\eta}$  
of closed subsets of $K$ with $K_0 = K$, $K_\eta = \emptyset$,  
$K_\gamma = \bigcap_{\alpha <\gamma} K_\alpha$ if $\gamma$ is a limit  
ordinal and $K_\alpha \supseteq K_\beta$ if $\alpha <\beta$,  
such that $F\big|_{  
K_\alpha\setminus K_{\alpha+1}}$ is continuous for all $\alpha$.  
\iitem{d)} There exists a sequence  $(K_n)$ of closed subsets of $K$  
with $K_n\subseteq K_{n+1}$ for all $n$ such that $K=\bigcup_n K_n$  
and $F\big|_{K_n}$ is continuous for all $n$.   
\smallskip\par  

\demo Remark 4.10.  
Property (c) suggests the following index for $PS(K)$:  
$$I(F) = \inf \bigl\{ \eta < w_1 : \exists\ (K_\alpha)_{\alpha\le \eta}  
\ \hbox{ satisfying (c)}\bigr\}\ .$$  

\noindent {\it Proof of 4.9.}  
d) $\Rightarrow$ a): Let $(K_n)$ be as in d) and for $n\in\IN$  let  
$f_n \in C(K_n)$ be given by $f_n = F\big|_{K_n}$. By the  
Tietze extension theorem there exists an extension of $f_n$,   
$\widetilde f_n\in C(K)$, with  
$\Vert \widetilde f_n\Vert_\infty \le \Vert F\Vert_\infty$. Clearly  
$(\widetilde f_n)$ is pointwise stabilizing and has limit $F$.  

a) $\Rightarrow$ b): For $n\in\IN$ set  
$$L_n = \bigl\{ k\in L : f_m (k) = F(k)\ \hbox{ for }\ m\ge n\bigr\}$$  
where $(f_n) \subseteq C(K)$, $\Vert f_n\Vert  \le  \Vert F\Vert$   
and $(f_n)$ is pointwise stabilizing with limit $F$. Let $G=\bigcup_n  
\hbox{int} (L_n)$. Thus $G$ is open in $L$.   
Also by the Baire Category theorem,   
$G$ is dense in $L$.  

b) $\Rightarrow$ c): Let $K_0 =K$ and let $K_1 = \sim G_0$ where $G_0$  
is a dense open subset of $K$   
and $F$ is continuous on $G_0$. Now if $K_\alpha$ is defined  
choose $G_\alpha$, a dense open subset of    
$K_\alpha$, so that $F\big|_{K_\alpha}$ is  
continuous on $G_\alpha$ and set $K_{\alpha +1} = K_\alpha \setminus  
G_\alpha$. At limit ordinals $\gamma$, set $K_\gamma =\bigcap_{\alpha<\gamma}  
K_\alpha$. Since $K$  is a separable metric space,  
$K_\eta = \emptyset$ for some $\eta <w_1$.  

c) $\Rightarrow$ d): Let $(K_\alpha)_{\alpha \le\eta}$ be as in c).  
Let $\E_n\downarrow 0$ and for each $n$ set $K_{\alpha,n} = \{ k\in K_\alpha:  
d(k,K_{\alpha +1}) \ge \E_n\}$ where $d$ is the metric on $K$. Let  
$K_n = \bigcup_{\alpha <\eta} K_{\alpha,n}$. We note that $K_n$ is closed.  
Indeed let $(k_i) \subseteq K_n$ converge to $k$. Then there exists  
$\alpha <\eta$ so that $k\in K_\alpha$ but $k\notin K_{\alpha+1}$. We  
claim that $k_i\in K_{\alpha,n}$ for sufficiently large $i$ and thus  
$k\in K_{\alpha,n}$ since $K_{\alpha,n}$ is closed. To see this note first  
that if $k_i \notin K_\alpha$, then $d(k_i,k) \ge \E_n$.  Thus for large  
$i$, $k_i\in K_\alpha$ and (since $k\notin K_{\alpha+1}$) $k_i\notin  
K_{\alpha +1}$. Hence $k_i\in  K_{\alpha,n}$ for large $i$ (since  
the $K_{\alpha,n}$'s are  disjoint in $n$).   

Finally $F\big|_{K_n}$ is continuous, for if $(k_i) \subseteq K_n$  
and $(k_i)$ converges to $k\in K_{\alpha,n}$, then by the above argument  
$k_i \in K_{\alpha,n}$ for large $i$ and $F\big|_{K_{\alpha,n}}$ is  
continuous.~\qed  
\medskip  

We end this section with an improvement of Proposition~4.8 in a special case.  

\proclaim Proposition 4.11.  
Let $K$ be a compact metric space and let $F$ be a simple  
Baire-1 function on $K$. Then there exists $(f_n)\subseteq C(K)$  
converging pointwise to $F$ such that $\lb (f_n)\rb$ embeds into $C(L)$  
for some countable compact space $L$.\par  

\proof   
First we consider the case where $K$ is totally disconnected.  
Choose $\E_0>0$ so that if $F(k_1) \ne F(k_2)$, then $|F(k_1) - F(k_2)|  
>\E_0$. Let $K_\alpha = K_\alpha (F,\E_0)$ for $\alpha \le \eta$  
with $K_\eta = \emptyset$. By our choice of $\E_0$, $F\big|_{K_\alpha\setminus  
K_{\alpha+1}}$ is continuous (with respect to $K_\alpha$) for all  
$\alpha <\eta$.  

Choose a countable partition $(D_j)$ of $K$ into closed sets with  
the following properties.  
\smallskip  
\iitem{a)} $\diam D_j \to 0$  
\iitem{b)} for each $j$, $D_j$ is a relatively clopen subset of $K_\alpha  
\setminus K_{\alpha+1}$ for some $\alpha <\eta$ such that $F\big|_{D_j}$  
is constant.  
\smallskip  
\noindent This can be done as follows. For each $\alpha$ choose a finite  
partition of relatively clopen subsets of $K_\alpha\setminus K_{\alpha+1}$  
such that $F$ is constant on each set of the partition. Each such set  
is relatively open in $K_\alpha$ and thus may be in turn  
partitioned into a countable number of relatively clopen subsets of   
$K_\alpha$.  List all the sets thus obtained for all $\alpha <\eta$  
as $(C_i)_{i=1}^\infty$. Each $C_i$ is closed in $K$ and thus may in  
turn be partitioned into a finite number of closed subsets of diameter  
not exceeding $1/i$. We list all these sets as $(D_j)_{j=1}^\infty$.  

Let $L= K/\{ D_j\}$ be the quotient space of $K$. Since each $D_j$  
is closed and $\diam D_j\to 0$, $L$ is compact metric. For $n\in \IN$  
choose $\hat f_n \in C(L)$ with $\Vert \hat f_n\Vert_\infty \le  
\Vert F\Vert_\infty$ and $\hat f_n(D_j)$ equal to the constant  
value of $F\big|_{D_j}$ for $j\le n$. Let $\phi : K \to L$ denote  
the quotient map and let $f_n = \hat f_n \circ \phi$. Clearly  
$f_n \in C(K)$, $\Vert f_n\Vert \le \Vert F\Vert$ and $(f_n)$  
converges pointwise to $F$. Also $[(f_n)]$ is isometric  
to $\lb (\hat f_n)\rb \subseteq C(L)$.  

For the general case let $\phi :\Delta\to K$ be a continuous surjection and  
let $F$ be a simple Baire-$1$ function on $K$. By the first part of the proof  
there exist $(f_n)\subseteq C(\Delta)$ converging pointwise to $F\circ\phi$  
and a countable compact metric space $L$ such that $\lb (f_n)\rb    
\hookrightarrow C(L)$. Let $(g_n)$ be a bounded sequence in $C(K)$ converging  
pointwise to $F$. By Lemma~3.3 there exist convex block subsequences  
$(h_n)$ and $(d_n) $ of $(g_n)$ and $(f_n)$, respectively, such that  
$\sum \| g_n\circ \phi - d_n\| <\infty$. Thus $\lb (g_n)\rb    
\cong \lb (g_n\circ \phi)\rb    
\hookrightarrow C(L)$.~\qed  

\demo  Question 4.12.  
Does Proposition~4.11 remain true   
if we only assume  
$F\in PS (K)$ or even $F\in DSC(K)$?  
Note that if $F$ satisfies the conclusion of 4.11, $F$ strictly  
governs the class of $c_0$-saturated spaces, while it is not clear  
that $DSC$ functions have this property.\par  

\beginsection{5. The Baire-1 Solar System.}  

In this section we shall examine the relationships between the various  
classes of Baire-1 functions which we have defined. We begin with a  
result which follows easily from the Banach space theory --- that   
developed above and some examples presented  in later sections.  

\proclaim Proposition 5.1.  
Let $K$ be an uncountable compact metric space. Then  
$$C(K) \subsetneqq DBSC (K) \subsetneqq B_{1/4} (K) \subsetneqq B_{1/2}(K)  
\subsetneqq B_1 (K)\ .\leqno(5.1)$$\par  

\proof Since $C(K)$ and $C(K')$ are isomorphic whenever $K$ and $K'$  
are both uncountable compact metric spaces \rf{29}, it suffices to  
separately consider each of the inclusions in (5.1).   
Thus if we show $C(K')\ne DBSC(K')$ for some uncountable compact metric  
space $K'$, then $C(K) \ne DBSC(K)$ as well. Indeed if $j: C(K)\to C(K')$  
is an onto isomorphism, then $\widetilde\jmath = j^{**}|_{B_1(K)}: B_1 (K)  
\to B_1(K')$. is an onto isomorphism satisfying $\widetilde\jmath (DBSC(K))  
= DBSC(K')$, $\widetilde\jmath (B_{1/4}(K)) = B_{1/4} (K')$ and  
$\widetilde\jmath (B_{1/2}(K)) = B_{1/2}(K')$.  

For the first inclusion, $C(K) \subsetneqq DBSC(K)$, let $X = c_0$.  
Then $K= (Ba (X^*),w^*)$ is uncountable compact metric and, as is well known,  
$X^{**} \subseteq DBSC(K)$. In particular if $F\in X^{**}\setminus X$,  
then $F\in DBSC(K)\setminus C(K)$.  

The fact that $B_{1/4} (K)\supsetneqq DBSC(K)$ follows from   
Theorem~A(b) and our example in \S6 where we produce a nonreflexive  
separable Banach space $X$ not containing $c_0$ such that $X^{**}  
\subseteq B_{1/4} (K)$, where  $K= Ba (X^*)$.  

For the next inclusion let $X= J$, the James space. $J$ is not  
reflexive and has no spreading model isomorphic to $c_0$ or $\ell_1$  
\rf{1}. Thus if $K = (Ba(J^*), w^*)$, then $X^{**}\setminus X \subseteq  
B_{1/2} (K)\setminus B_{1/4} (K)$ by virtue of Theorem~B.  

For the last inclusion let $Y$ be the quasi-reflexive space of order~1  
(see the proof of Proposition~6.3) whose dual is $J(e_i)$, where  
$(e_i)$ is the unit vector basis of Tsirelson's space.  It is proved in  
\rf{32} that the only spreading models of $Y$  are isomorphic to $\ell_1$.  
Thus by Theorem~3.7, if $Y^{**} = Y\oplus \langle F\rangle$ and  
$K= Ba (Y^*)$, then $F\notin B_{1/2} (K)$. An alternative method  
would be to consider the quasi-reflexive spaces $Q_\gamma$ constructed  
in Proposition~2.8.~\qed  

\demo Remark 5.2.  
How does the class $DSC(K)$ relate to the classes in (5.1)? Of course we  
always have $DBSC(K) \subseteq DSC(K)$ and in fact for $K$ an  
uncountable compact metric space we have the following diagram.  
$$\matrix{&\cr \noalign{\vskip2.2truein}&\cr}$$  

\noindent   
Thus $DSC(K)$ is an asteroid in the Baire-1 solar system. Indeed our proof  
of Proposition~5.1 along with Theorem~4.1 yields  that  
$B_{1/4}(K) \setminus DSC(K) \ne \emptyset$, $B_{1/2} (K)\setminus  
\lb DSC(K)\cup B_{1/4}(K)\rb \ne \emptyset$ and     
$B_1 (K) \setminus \lb DSC(K) \cup  
B_{1/2}(K)\rb \ne\emptyset$. The fact that $DSC(K)\cap B_1 (K)\setminus  
B_{1/2} (K)$, $DSC (K)\cap B_{1/2} (K)\setminus B_{1/4}(K)$ and  
$DSC(K) \cap B_{1/4}(K)\setminus DBSC(K)$ are all nonempty follows  
from Proposition~5.3 below.  

We now turn to the case  where  $K$ is a {\it countable\/} compact  
metric space. In this setting we have, of course, $DSC(K) = B_1(K)$.  
However if $K$ is large enough, the classes in (5.1) are still distinct.  
Since every countable compact metric space is homeomorphic to some  
countable ordinal, given the order topology \rf{30}, we confine ourselves  
to this setting.  

\proclaim Proposition 5.3.  
\smallskip  
\iitem{a)} If $K = \omega^{\omega^2} +$, then $B_{1/4}(K)\setminus  DBSC(K)  
\ne\emptyset$.  
\iitem{b)} If $K= \omega^\omega +$, then $B_{1/2} (K)\setminus B_{1/4} (K)  
\ne\emptyset$.  
\iitem{c)} If $K= \omega^\omega +$, then $B_1 (K)\setminus B_{1/2}(K)  
\ne\emptyset$.  
\iitem{d)} If $K= \omega^+$, then $DBSC(K)\setminus C(K) \ne\emptyset$.  
\smallskip\par  

Before proving this proposition we need some terminology. Recall that an  
indicator function $\bone_A$ is Baire-1 iff $A$ is ambiguous (simultaneously  
$F_\sigma$ and $G_\delta$). Thus if $A\subseteq K$ where $K$ is   
countable compact metric, then $\bone_A\in B_1 (K)$. We begin with a   
discussion of such functions.  

Let $\delta$ be a countable compact ordinal space (in its order topology).  
Recursively we define $I_0=\emptyset$,   
$I_1 = \{ x\in\delta : x$ is an isolated point of  
$\delta\}$, and for $\alpha >1$, $I_\alpha = \{ x\in\delta\setminus  
\bigcup_{\beta <\alpha } I_\beta : x$ is an isolated point of $\delta\setminus  
\bigcup_{\beta <\alpha}I_\beta\}$. The $I_\alpha$'s are just the relative 
complements of the usual derived sets. 

Let us say an ordinal is {\it even\/} if it is of the form $\gamma + 2n$  
for some $n\in\IN$ where $\gamma=0$ or  
$\gamma $ is a  limit ordinal. Let $F_\delta = \bone_{  
A_\delta}$ where $A_\delta = \bigcup_{\alpha\, \hbox{\sevenrm even}}  
I_\alpha$. We have  
\smallskip  
\iitem{$^\circ 1$)} $\Vert F_{\omega^n+}\Vert_\infty =1$ and  
$|F_{\omega^n +} |_D =n$.  
\iitem{$^\circ 2$)} $|F_\delta |_D =\infty$ if $\delta\ge \omega^\omega +$.  
\smallskip  
\noindent $^\circ 1$) implies $^\circ 2$) trivially. To see $^\circ 1$),  
one first notes that $K_n (F_{\omega^n +},1) \ne \emptyset$. Indeed,  
$K_\alpha (F_\delta,1)$ is just the $\alpha^{th}$ derived set of $\delta$.  
Hence $|F_{\omega^n+}|_D\ge n$ by the proof of  
Lemma~2.4.  We leave the reverse  
inequality to the reader.  

\demo Definition.  
We say that a function $F:\omega^n + \to \IR$ is of {\it type $0$\/} if  
$F= n^{-1} F_{\omega^n+}$. The  domain of $F$, $\omega^n+$, is called  
a {\it space of type $0$\/}.  

Thus if $F$ is a function of type 0 with domain $\omega^n+$,   
$|F|_D =1$ and $\Vert F\Vert_\infty =n^{-1}$.    

More generally for $n\in\IN $ we have the   

\demo Definition.  
A class of real valued functions $\F_n$ defined on countable compact  
metric spaces is said to be of {\it type $n$\/} if  
\smallskip  
\iitem{a)} For $F\in \F_n$, $|F|_D\ge n$.  
\iitem{b)} For $F\in \F_n$, $F$ is the uniform limit of $(F_m)$ with  
$\sup_m |F_m|_D\le 1$.  
\iitem{c)} For each $\varep  >0$, there is an $F\in \F_n$ with  
$\Vert F\Vert_\infty < \varep$.  
\smallskip\par  

The domain of $F\in\F_n$ is called a {\it space of type $n$.}  

\proclaim Lemma 5.4.  
For $n\in\IN\cup \{ 0\}$ there exists a class $\F_n$ of functions  
of type $n$.\par  

\proof We have seen that $\F_0$ exists. Suppose $\F_n$ exists. To  
obtain functions $F\in\F_{n+1}$ we begin with a function $G\in \F_0$  
defined on a set $K$. Let $(t_i)_{i=1}^\infty$ be a list of the isolated  
points of $K$. We enlarge $K$ as follows. To each $t_i$ we adjoin  
a sequence of disjoint clopen sets $K_1^i,K_2^i,\ldots\ $ clustering  
only at $t_i$. Each of the $K_j^i$'s is a space of type~$n$ supporting  
a function $F_j^i$ of type~$n$ with $\Vert F_j^i\Vert_\infty \le  
(i+j+m)^{-1}$. Here $m\in\IN$ is arbitrary but fixed. $K_{n+1}$, the  
new space of type $n+1$, is this enlarged space. Set  
$$F(t) = \cases{ G(t)\ ,&\quad $t\in K$\cr  
\noalign{\vskip4pt}  
	F_j^i (t)\ ,&\quad $t\in K_{ij}$\ .\cr}$$  
Let $\F_{n+1}$ be the set of all such $F$'s thusly obtained. We must  
check that $\F_{n+1}$ satisfies a) and b) with $n$ replaced by $n+1$  
(~c) is immediate). b) holds since $F$ is the uniform limit of  
$(F_k)$ where  
$$F_k(t) = \cases{ G(t)\ ,&\quad $t\in K$\cr  
\noalign{\vskip4pt}  
F_j^i(t)\ ,&\quad $t\in K_j^i$  with $i+j \le k$\cr  
\noalign{\vskip4pt}  
0\ ,&\quad otherwise\cr}$$  
and each $F_k$ is the uniform limit of $(F_{k,n})_{n=1}^\infty$  
where $|F_{k,n}|_D \le 1$ for all $n$.  

Finally we check a). Let $(f_m)_0^\infty \subseteq   
C(K_{n+1})$, $f_0\equiv 0$,  
converge pointwise to $F$. Since $F\big|_K = G$   
and $|G|_D =1$, for $\varep >0$  
there exist $t_{i_0} \in K$ and $k\in\IN$ with   
$\sum_{i=0}^{k-1} |f_{i+1} (t_{i_0}) - f_i (t_{i_0}) | > 1-\varep$.  
Moreover by the nature of $G$ we may assume $|G(t_{i_0})|<\varep$.  
Since the $K_j^{i_0}$'s cluster at $t_{i_0}$ and each $f_i$ is   
continuous there exists $j_0\in \IN$ so that for $t\in K_{j_0}^{i_0}$,  
$\sum_{i=0}^{k-1} |f_{i+1} (t) - f_i(t)| > 1-\varep$ and $|f_k(t)| <\varep$.  
But on $K_{j_0}^{i_0}$, $(f_m)$ converges pointwise to $F_{j_0}^{i_0}$  
and $|F_{j_0}^{i_0}|_D\ge n$. Thus there exists $t\in K_{j_0}^{i_0}$ with  
$$|f_{k+1}(t)| + \sum_{i>k} |f_{i+1} (t) -f_i (t)| > n-\varep\ .$$  
It follows that   
$$\sum_{i=0}^\infty | f_{i+1} (t) - f_i (t)| > n+1 - 3\varep$$  
which proves a).~\qed  

\demo Remark 5.4.  
Our proof yields that the spaces of type-$n$ can be constructed within   
$\omega^{\omega \cdot(n+1)}+$.\par  

\noindent {\it Proof of Proposition 5.3.}  
a) Let $K= \omega^{\omega^2}+$ and choose (by Remark~5.4) a sequence  
$(K_n)_{n=0}^\infty$ of disjoint clopen subspaces of $K$ with  
$K_n$ of type-$n$. Let $F_n$ be a function of type-$n$ supported  
on $K_n$ with $\Vert F_n\Vert_\infty \to 0$ and let $F$ be the sum of  
the $F_n$'s. Clearly $|F|_D=\infty$ since $ |F_n|_D \ge n$. Yet $F$  
is the uniform limit of a sequence of  functions with $|\cdot|_D$ not  
exceeding~1.  

b) Let $(K_n)_{n\in\IN}$ be a sequence of disjoint clopen subspaces of  
type-0 of $\omega^\omega + = K$ such that $K_n$ supports a function  
$F_n$, which is a multiple of a function of type-0, with $\Vert F_n  
\Vert_\infty \le n^{-1}$ and $|F_n|_D \ge n$. Define   
$$F(t)  =\cases{ F_n(t) &if $t\in K_n$\cr  
\noalign{\vskip4pt}  
	0&otherwise.\cr}$$  
Clearly $F\in B_{1/2} (K)\setminus B_{1/4}(K)$.  

c) The type-0 function $F_{\omega^\omega +}$ is not Baire-1/2.  

d) $F_{\omega^+}$ is $DBSC$.~\qed  
\medskip  

It is easy to check that the results of Proposition~5.3 are best possible.  

\beginsection{6. A Characterization of $B_{1/4}(K)$ and an Example.}  

In this section we give an example which shows that functions of class  
Baire-1/4 need not govern $\{ c_0\}$. Thus Theorem~B(b) is best possible.  
Before giving the example we give a sufficient (and necessary) criterion   
for a function to be Baire-1/4.  

\proclaim Theorem 6.1.  
Let $K$ be a compact metric space and let $F\in B_1(K)$.   
Then $F\in B_{1/4}(K)$  
iff there exists a $C<\infty$ such that for all $\varep >0$ there exists  
a sequence  
$(S_n)_{n=0}^\infty \subseteq C(K)$, $S_0\equiv 0$, with $S_n(k)\to F(k)$  
for all $k\in K$ and such that for all subsequences $(n_i)$ of   
$\{0\}\,\cup\,\IN$  
and $k\in K$,  
$$\sum_{j\in B((n_i),k)} |S_{n_{j+1}} (k) - S_{n_j} (k)| \le C\ .  
\leqno(6.1)$$  
Here $B((n_i),k) = \{ j:|S_{n_{j+1}} (k) - S_{n_j} (k)| \ge \varep\}$.\par  

\proof First  assume $F\in B_{1/4} (K)$, let $\varep >0$ and let $\varep_n  
\downarrow 0$. By the proof of Theorem~B(b) there exists $(f_n)_{n=0}^\infty  
\subseteq C(K)$, $f_0\equiv 0$, converging pointwise to $F$ with the   
following property. For each $m\in\IN$, there exists $(h_j^m)_{j=0}^\infty  
\subseteq C(K)$ with $h_0^m\equiv 0$ and   
$$\sum_{j=0}^\infty |h_{j+1}^m (k) - h_j^m(k)| \le M \equiv 2|F|_{1/4}\ ,  
\ \hbox{ for }\ k\in K\ .\leqno(6.2)$$  
Furthermore $\Vert h_j^m - f_j\Vert_\infty \le \varep_m$ for $j\ge m$.  

Let $\varep >0$ and fix $m$ with $4\varep_m <\varep$.   
Let $(S_n)_{n=0}^\infty = (0,f_m,f_{m+1},\ldots)$, and let  
$(n_i)$ be a subsequence of $\{ 0\} \cup \IN$ and let $k\in K$  
be fixed. Then  
$$\sum_{j\in B((n_i),k)} |S_{n_{j+1}} (k) - S_{n_j} (k)| \le  
\sum_{j=0}^\infty | h_{j+1}^m (k) - h_j^m(k)| + 2\varep_m\# B((n_i),k)\ .  
\leqno(6.3)$$  
Since $|f_p (k) - f_q(k)| \ge \varep$ implies for $p>q\ge m$ or $q=0$  
that $|h_p^m (k) - h_q^m (k) | \ge\varep -2\varep_m> \varep/2$,   
(6.2) yields that $\# B  
((n_i),k)\le 2M/\varep$. Thus (6.3) yields (6.1) with $C= 2M = 4|F|_{1/4}$.  

For the converse, let $C>\varep >0$ and let $(S_n)_0^\infty \subseteq  
C(K)$, $S_0 \equiv 0$, converge pointwise to $F$ and satisfy (6.1)  
for any subsequence $(n_i)$ of $\{ 0,1,2,\ldots\}$ and any $k\in K$.   
For $k\in K$ we linearly extend the sequence $(S_n(k))_{n=0}^\infty$  
to $(S_r(k))_{r\ge0}$. Precisely, if $r=\lambda n+ (1-\lambda) (n+1)$  
we set $S_r (k) = \lambda S_n (k) + (1-\lambda) S_{n+1}(k)$. Since  
the $S_n$'s are continuous, $S_r \in C(K)$ as well. Furthermore, if  
$0\le r_1<r_2<r_3<\cdots\ $,  $k\in K$ and $B= B((r_i),k) = \{ j:  
|S_{r_{j+1}} (k) - S_{r_j} (k)| \ge \varep\}$, then   
$$\sum_{j\in B} |S_{r_{j+1}} (k) - S_{r_j} (k)| \le 3C\ .\leqno(6.4)$$  
Indeed if $J_n = \{ j\in B: n\le r_j < r_{j+1} \le n+1\} \ne \emptyset$,  
then $\varep \le \sum_{j\in J_n} |S_{r_{j+1}} (k) - S_{r_j} (k) | \le  
| S_{n+1} (k) - S_n (k)|$.  
If $j\in B\setminus \bigcup_n J_n$, there exists integers $\ell_j$  
and $m_j$ with $\ell_j - 1\le r_j < \ell_j \le m_j < r_{j+1} \le  
m_{j+1}$. Thus by linearity for some choice of $p_j \in \{ \ell_j-1,  
\ell_j\}$ and $q_j \in \{ m_j,m_j+1\}$ we have $\varep \le |S_{r_{j+1}}  
(k) - S_{r_j} (k) | \le |S_{q_j} (k) - S_{p_j} (k)|$. Thus  
$$\eqalign{&\sum_{j\in B} |S_{r_{j+1}}(k) - S_{r_j} (k) | \le   
	\sum_{\{ n: J_n\ne \emptyset\}} |S_{n+1}(k) -S_n(k)| \cr  
&\qquad +\sum_{2j\in B\setminus \bigcup_n J_n} | S_{q_{2j}} (k) -  
S_{p_{2j}} (k)| + \sum_{2j+1\in B\setminus\bigcup_nJ_n}  
|S_{q_{2j+1}} (k) - S_{p_{2j+1}} (k)| \le 3C\ .\cr}$$  

We shall construct a sequence $(f_n)_{n=0}^\infty \subseteq C(K)$,  
$f_0\equiv 0$, such that for $k\in K$,  
$$\leqalignno{& \sum_{n=0}^\infty |f_{n+1} (k) - f_n(k)| \le 4C\ \hbox{ and }  
	&(6.5)\cr  
&\hbox{ if $H$ is the pointwise limit of $(f_n)$ then $\Vert H-F\Vert_\infty  
\le 5\varep$\ .}&(6.6)\cr}$$  
This will complete the proof.  
\medskip  

Each $f_n$ shall be an average of functions $S_t$ where $t:K\to\lb 0,\infty)$  
is continuous and $S_t (k) \equiv S_{t(k)} (k)$ for $k\in K$.  
Let $f_0 = S_0\equiv 0$. Let $\alpha_1^1 : \lb 0,\infty) \to \lb 0,1\rb $ be  
identically 0 on $\lb 0,\varep\rb $, identically 1 on 
$\lb 3\varep /2,\infty)$  
and linear on $\lb\varep,3\varep /2\rb$. Let $\alpha_2^1 : \lb 0,\infty) \to  
\lb 0,1\rb$ be identically 0 on $\lb 0,3\varep /2\rb$, identically 1 on  
$\lb 2\varep,\infty)$ and linear on $\lb 3\varep /2,\varep\rb$.   
For $i=1,2$ let $t_i (k) = \alpha_i^1 (|S_1(k)|)$. Let $f_1 = 2^{-1}  
(S_{t_1} + S_{t_2})$. We next define continuous functions $t_{i,j}$ for  
$i=1,2$ and $j= 1,2,3,4$ by $t_{i,j} (k) = t_i(k)  + \alpha_j^2  
(|S_2(k)- S_{t_i}(k)| ) (2-t_i(k))$. Here $\alpha_j^2 : \lb 0,\infty)\to  
\lb 0,1\rb$ is identically 0 on $\lb 0,(4+j-1)\varep /4\rb$, identically 1 on  
$\lb (4+j)\varep /4,\infty)$ and linear on     
$\lb (4+j-1)\varep /4,(4+j)\varep/4\rb$.  
Set $f_2 = 8^{-1}\sum_{i=1}^2 \sum_{j=1}^4 S_{t_{i,j}}$.  

In general if  $f_n =2^{-1} 2^{-2}\cdots 2^{-n}\sum S_{t_{i_1,\ldots,i_n}}$,  
where the  indices of summation range over  $\{ (i_1,\ldots, i_n) :  
1\le i_j \le 2^j\}$. We define $t_{i_1,\ldots, i_n}$ for $1\le i_{n+1}\le  
2^{n+1}$ by   
$$\eqalign{& t_{i_1,\ldots,i_{n+1}}(k) = t_{i_1,\ldots,i_n} (k)\cr  
&\qquad + \alpha_{i_{n+1}}^{n+1} \bigl( |S_{n+1}(k) - S_{t_{i_1,\ldots,i_n}}  
	(k)|\bigr) \bigl( n+1 -t_{i_1,\ldots,i_n}(k)\bigr)\ .\cr}$$  
The functions $\alpha_j^{n+1}$ for $1\le j\le 2^{n+1}$ are defined as  
before to be  identically 0 on $\lb 0,\varep + (j-1) \varep 2^{-n-1}\rb $,  
identically 1 on $\lb \varep +j\varep 2^{-n-1}\rb$ and linear on  
$\lb \varep + (j-1)\varep 2^{-n-1},\varep +j\varep 2^{-n-1}\rb$.  

The  point of the construction is this. For $k\in K$ and $(i_1,\ldots,i_n)$  
fixed, $|S_{t_{i_1,\ldots,i_{n+1}}}(k) - S_{t_{i_1,\ldots, i_n}}(k)|$  
is either 0 or a number exceeding $\varep$ for all but perhaps one choice  
of $i_{n+1}$. $\lb$This is because the nonconstant     
parts of the $\alpha_j^{n+1}$'s  
are disjointly supported.$\rb$     
Also except for at most one value of $i_{n+1}$,  
if $|S_{t_{i_1,\ldots, i_{n+1}}}(k) - S_{t_{i_1,\ldots, i_n}}(k)| \ge  
\varep$ then $S_{t_{i_1,\ldots,i_{n+1}}}(k) = S_{n+1} (k)$.   

We next check (6.5). Fix  $k\in K$ and $m\in\IN$. A simple calculation using  
the triangle inequality shows that  
$$\sum_{n=0}^m |f_{n+1} (k) - f_n(k)| \le \hbox{AVE}\sum_{n=0}^m  
| S_{t_{i_1,\ldots, i_{n+1}}} (k) - S_{t_{i_1,\ldots, i_n}}(k)|  
\leqno(6.7)$$  
where the average is taken over $\{ (i_1,\ldots, i_m) : 1\le i_j\le 2^j$  
for all $j\}$. If we fix $(i_1,\ldots, i_m)$ and let   
$$B = \bigl\{ n\le m\mid | S_{t_{i_1,\ldots, i_{n+1}}}(k)  
- S_{t_{i_1,\ldots, i_n}}(k) | \ge \varep\bigr\}$$  
then  
$$\sum_{n\in B} | S_{t_{i_1,\ldots, i_{n+1}}}(k) - S_{t_{i_1,\ldots,i_n}}  
(k)| \le 3C$$   
by (6.4).  

Now for $1\le n\le m$ fixed, the percentage of terms in the ``AVE'' of (6.7)  
for which $0<|S_{t_{i_1,\ldots, i_{n+1}}}(k) - S_{t_{i_1,\ldots,i_n}}(k)|  
<\varep$ is at most $2^{-n-1}$. It follows that  
$$\hbox{AVE }\sum_{n=0}^m |S_{t_{i_1,\ldots, i_{n+1}}} (k)  
- S_{t_{i_1,\ldots, i_n}}(k) | \le 3C + 2^{-1} \varep  
+\cdots + 2^{-m-1} \varep$$  
and (6.5) follows from this since $\varep <C$.  

(6.5) implies $(f_n)$ is pointwise convergent to some function $H$. For fixed  
$k\in K$ choose $m\in \IN$ so that $2^{-m} C<\varep$, $|S_m (k) - F(k)| <  
\varep$ and $|f_m (k) - H(k)| <\varep$. We claim that $|f_m(k) - S_m(k)| <3  
\varep$, which proves (6.6). Indeed  
$$f_m (k) = \hbox{AVE } S_{t_{i_1,\ldots, i_m}}(k) \ \hbox{ and }\   
C\ge |S_{t_{i_1,\ldots, i_m}}(k) - S_m (k)| \ge 2\varep$$  
for at most $2^{-m}\# \{ (i_1,\ldots, i_m) : i_j \le 2^j\}$ choices  
of $(i_1,\ldots, i_m)$.  
Thus $|f_m(k) - S_m(k)| \le 2\varep + 2^{-m} C < 3\varep$.~\qed  

\demo Remark 6.2.  
Let $F\in B_1 (K)$. Our proof shows that $F\in B_{1/4}(K)$ iff there exists  
$C<\infty$ and $(S_n)_{n=0}^\infty \subseteq C(K)$, $S_0 \equiv 0$,  
converging pointwise to $F$ such that for all $\varep >0$ there exists  
$m\in\IN$ such that if $(n_i)$ is any subsequence of $\{ 0,m,m+1,\ldots\}$  
then (6.1) holds.\par  

\proclaim Proposition 6.3.  
There exists a compact metric space $K$ and $F\in B_{1/4}(K)$ which does   
not govern $\{ c_0\}$.\par  

\proof Let $(e_i)$ be the unit vector basis of the Tsirelson  
space $T$ constructed in \rf{17} (see also \rf{11}) and let  
$X= J(e_i)$ be its ``Jamesification'' as described in \rf{6}.   
For completeness we recall the definition of $X$. Let $c_{oo}$ be the   
linear space of all finitely supported functions $x: \IN\to\IR$ and for  
$n\in\IN$ define $S_n : c_{oo} \to \IR$ by $S_n(x) = \sum_{i=1}^n x(i)$.  
Let $S_0 \equiv 0$. For $x\in c_{oo}$ let  
$$\Vert x\Vert = \sup \biggl\{ \Big\Vert \sum_{i=1}^m (S_{n_i} - S_{p_i-1})  
(x) e_{p_i} \Big\Vert_T \ \Big|\    
1\le p_1 \le n_1 < p_2 \le n_2 < \cdots < p_m \le n_m\biggr\}\ .$$  
Let $X$ be the completion of $(c_{oo},\Vert\cdot\Vert)$.  

As shown in \rf{6}, the unit vectors $(u_i)$ form a boundedly complete  
normalized basis for $X$. Thus $X= Y^*$ (where $Y = \lb (u_i^*)\rb    
\subseteq X^*$). Furthermore it was shown that $Y$ is quasi-reflexive  
and $Y^{**}$ has a basis given by $\{ S,u_1^*, u_2^*,\ldots \}$, where  
$$S\Bigl(\sum a_iu_i\Bigr) = \sum_1^\infty a_i\ .$$  
Of course $(u_i^*)$ are the biorthogonal functionals to $(u_i)$ and  
$S$ is the weak* limit in $Y^{**}$ of $(S_n)$.  

Let $K= Ba (X) = Ba (Y^*)$ in the weak* topology (of $Y^*$). Since $Y$  
does not contain $c_0$, our example will be complete if we can prove  
that $S\in B_{1/4}(K)$. By Theorem~6.1 it suffices to prove that if  
$\varep >0$ then for $m\in \IN$ with $m > 2/\varep$, if $x\in Ba(X)$  
and $(n_i)$ is a subsequence of $\{ m,m+1,m+2,\ldots\}$, then  
$$\sum_{j\in B} | S_{n_{j+1}} (x) - S_{n_j} (x) | \le 2$$  
where  
$$B = \bigl\{ j: |S_{n_{j+1}} (x) - S_{n_j} (x) | \ge \varep\bigr\}\ .$$  

We first note that $\# B<m$. Indeed if $\# B\ge m$, then by the   
properties of $T$,   
$$\eqalign{  
1&\ge \Vert x\Vert \ge \Big\Vert \sum_{j\in B} \bigl( S_{n_{j+1}} (x) -  
S_{n_j} (x)\bigr) e_{n_j}\Big\Vert_T \cr  
\noalign{\vskip6pt}  
&\ge 2^{-1} m\varep\ ,\cr}$$  
a contradiction. The last inequality is due to the fact that  
$\Vert \sum_A a_i e_i\Vert_T \ge 2^{-1} \sum_A | a_i|$ provided  
$\min A \le \# A$.  

Thus $m\le \min B \le \# B$ and so  
$$\eqalign{\sum_{j\in B} |S_{n_{j+1}} (x) - S_{n_j} (x) |   
	&\le 2\Big\Vert \sum_{j\in B} \bigl( S_{n_{j+1}} (x)  
	- S_{n_j} (x)\bigr) e_{n_j} \Big\Vert \cr  
\noalign{\vskip6pt}  
&\le 2\Vert x\Vert \le 2\ .\cr}  
\eqno\blackbox$$  

\demo Remark 6.4.  
Our proof of Proposition 6.3 shows that there exists a quasi-reflexive  
(of order one) Banach space $Y$ such that if $K= Ba (Y^*)$ then  
$Y^{**} \setminus Y\subseteq B_{1/4} (K)$. In particular, it follows  
that there exists an $F\in B_{1/4} (K)\setminus C(K)$ which  
strictly governs the class of quasi-reflexive Banach spaces.\par

\beginsection{7. Some Bad Baire-1/2 Functions.}  

In this section we show that functions   of class Baire-1/2 need  
not be that nice.  

\proclaim Proposition 7.1.  
There exists a compact metric space $K$ and $F\in B_{1/2}(K)$ which   
governs $\{ \ell_1\}$.\par  

\demo Remark 7.2.  
The first example of an $F\in B_1 (K)$ which governs $\{ \ell_1\}$ was  
due to Bourgain \rf{9,10}. His ingenious construction forms the  
motivation behind our next example (Proposition~7.3). Another example  
of such an $F$ appears in \rf{2}. While the example of \rf{2} can be  
shown to be Baire-1/2, we prefer to present a very slight modification.  

\proof Let $(e_n)$ be the unit vector basis of a Lorentz sequence space  
$d_{w,1}$ (see {\it e.g.}, \rf{27}). Let $J(e_i)$ be the Jamesification  
of $(e_n)$ (see \rf{6}) and let $(u_i)$ be the unit vector basis of  
$J(e_i)$. Thus  
$$\Big\Vert \sum_{i=1}^k a_i u_i\Big\Vert = \sup  
\left\{\Big\Vert \sum_{i=1}^p \biggl( \sum_{j=n_i}^{m_i}    
a_j\biggr) e_i\Big\Vert_{d_{w,1}}\ \Big|\  
{1\le n_1\le m_1<n_2\le m_2 <\cdots < n_p\le m_p}\right\}  
\ .$$  
$(u_i)$ is a normalized spreading basis for $J(e_i)$ which is not  
equivalent to the unit vector basis of $\ell_1$ and thus by \rf{36},  
$(u_i)$ is weak Cauchy. Furthermore by standard block basis arguments  
one can show that $J(e_i)$ is hereditarily $\ell_1$. Also if $F$ is defined  
by $u_i \to F$ weak* then $F\in B_{1/2}(K)$ where   
$K= Ba (J(e_i)^*)$. But this is immediate by Theorem~B(a) since   
$(u_i)$, being its own spreading model, does not have $\ell_1$ as  
a spreading model. The fact that $F$ governs $\ell_1$ follows from  
Lemma~3.3. Indeed if $(f_n)$ is a bounded sequence in $C(K)$ converging  
pointwise to $F$, then some convex block subsequence of $(f_n)$ is a  
basic sequence equivalent to some convex block subsequence of  
$(u_i)$. Since $\lb (u_i)\rb$ is hereditarily $\ell_1$, $\ell_1\hookrightarrow  
\lb (f_n)\rb$.~\qed  

\proclaim Proposition 7.3.  
There exists a compact metric space $K$ and $F\in B_{1/2} (K)$ such that  
$F$ does not govern $\{ \ell_1\}$ yet $F$ strictly   
governs $\{ X : X$ is separable  
and $X^*$  is not separable$\}$.\par  

\demo Remark 7.4.  
In \rf{33} a function $F\in B_1 (K)\setminus B_{1/2}(K)$ was constructed  
satisfying the conclusion of Proposition~7.3. The construction we now  
present will be a modification of that example.\par  

\noindent {\it Proof of Proposition 7.3.}  
We begin by defining a Banach space $Y$. (The space $Y$ was first defined  
in \rf{34}) Let $\D = \{ \phi\} \cup\bigcup_n \{ 0,1\}^n$ be the   
dyadic tree with its natural order (see Remark~4.2) and let   
$(K_\alpha)_{\alpha\in\D}$ be the natural clopen base for the Cantor set  
$\Delta$. For $f\in C(K_\alpha)$ we let $\widetilde f \in C(\Delta)$  
be given by $\widetilde f (t) = f(t)$ for $t\in K_\alpha$ and   
$\widetilde f(t) = 0$ otherwise. Let  
$$\eqalign{  
&Y= \biggl\{ (f_\alpha)_{\alpha\in\D}\mid f_\alpha \in C(K_\alpha)  
\ \hbox{ for all }\ \alpha \in \D\ \hbox{ and }\cr  
\noalign{\vskip6pt}  
&\Vert (f_\alpha)\Vert_Y \equiv \sup \biggl\{ \biggl( \sum_{k=1}^\ell  
\Big\Vert \sum_{\alpha\in S_k} \widetilde f_\alpha\Big\Vert_\infty^2\biggr)  
^{1/2} :  
(S_k)_{k=1}^\ell\ \hbox{ are disjoint segments in }\ \D\biggr\}  
< \infty \biggr\}\ .\cr}$$  
$Y$ is a Banach space under the given norm.  

We shall construct a weak Cauchy sequence $(g_n)\subseteq Y$ with weak*  
limit $F$ such that   
$$\leqalignno{  
&\qquad \ell_1 \not\hookrightarrow \big\lb (g_n)\big\rb\ ,&(7.1)\cr  
\noalign{\vskip6pt}  
&\left\{ \eqalign{  
	&\big\lb (h_n)\big\rb^*\   
	\hbox{ is nonseparable for every convex block}\cr  
	&\qquad \hbox{subsequence $(h_n)$ of $(g_n)$ and} \cr}  
	\right.&(7.2)\cr  
\noalign{\vskip6pt}  
&\left\{ \eqalign{  
	&\hbox{there exists a weak* closed set $K\subseteq Ba (Y^*)$   
		such that}\cr  
	&\qquad K\ \hbox{ norms $\lb (g_n)\rb$ and }\   
	F\big|_K\in B_{1/2}(K)\ .\cr}  
\right.&(7.3)\cr}$$  
The proposition follows immediately from (7.1)--(7.3). Indeed to see that  
$F$ governs $\{ X: X$ is separable and $X^*$ is nonseparable$\}$, let  
$(f_n)$ be a bounded sequence in $C(K)$ converging pointwise to $F$. By  
Lemma~3.3 there exist convex block subsequences $(d_n)$ and $(h_n)$ of  
$(f_n)$ and $(g_n)$, respectively, such that $\Vert d_n-h_n\Vert_{C(K)} \to	
0$. Since $\lb (h_n)\rb^*$ is nonseparable, so is $\lb (d_n)\rb^*$.  

Our construction of $(g_n)$ depends upon the following (which in turn  
follows from our discussion of functions of type-0 in \S5):  
for $n\in \IN$ there exists $F_n \in B_1(\Delta)$ such that  
$$\Vert F_n\Vert_\infty =1\ \hbox{ and }\leqno(7.4)$$  
$$\left\{ \eqalign{& |F_n|_D = n\ .  
\hbox{ Moreover if $(h_i) \subseteq C(\Delta)$ converges  pointwise to $F_n$  
	then}\cr  
&\hbox{there exists $k\in\Delta$, integers $\ell_1 < \ell_2 <\cdots  
\ell_{n+1}$ and $\varep_i = \pm 1$ $(1\le i\le n)$}\cr  
&\hbox{such that  
$\sum\limits_{i=1}^n \varep_i (h_{\ell_{i+1}}- h_{\ell_i}) (k) > n-1$.}\cr}  
\right.  
\leqno(7.5)$$  

Actually our $F_n$'s are indicator functions whose domains  
are countable compact metric spaces $K$. Of course one can embed $K$ into  
$\Delta$ and the corresponding extended indicator functions have  
the desired properties (7.4) and (7.5).  

We use ``$<_L$'' for the natural linear order on $\D$.   
Thus $\phi < 0 < 1 < 00 < 01 < 10 < 11 < 000 < \cdots\ $.  
For each $\alpha  
\in \D$ choose $n_\alpha \in \IN$ and $c_\alpha \in \IR^+$ satisfying the  
following seven properties:  
\smallskip  
\iitem{i)} $\sum_{\beta \in \D} c_\beta \le 1$.  
\iitem{ii)} $c_\alpha^{-1} n_\alpha^{-1} \sum_{\beta <\alpha} n_\beta  
	< 1/10$.  
\iitem{iii)} $2c_\alpha^{-1} \sum_{\beta >\alpha} c_\beta < 1/10$.  
\iitem{iv)} $1-n_\alpha^{-1} > 9/10$.  
\iitem{v)} $2c_{\alpha_0} c_\alpha^{-1} < 1/10$ if $\alpha <_L \alpha_0$.  
\iitem{vi)} $c_{\alpha_0} c_\alpha^{-1} n_{\alpha_0} n_\alpha^{-1}  
< 1/10$ if $\alpha_0 <_L \alpha$.  
\iitem{vii)} $\sum_{\beta \in\D} b_\beta^2 < \infty$ where  
$b_\beta = \sum_{\gamma \ge_{_L}\beta} c_\gamma$.  
\smallskip  

\noindent Of course we could trim this list somewhat, but we prefer  
to list the properties in the form in which they are used. The   
$c_\alpha$'s and $n_\alpha$'s can be chosen as follows. Let   
$\{ \alpha_1,\alpha_2,\alpha_3,\ldots\}$ be a listing of $\D$  
in the linear order. Let $c_{\alpha_j} = (22)^{-j}$. It is quickly  
checked that properties i), iii), v) and vii) hold. We then choose  
$n_{\alpha_j}$ inductively to be an increasing sequence of positive  
integers with $n_{\alpha_1} = 11$ (so that iv) holds). If $n_{\alpha_j}$  
is picked, choose $n_{\alpha_{j+1}}$ to satisfy ii) and vi) for  
$\alpha = \alpha_{j+1}$.   
For each $\alpha\in\D$, let $F_{n_\alpha} \in B_1 (K_\alpha)$ satisfy  
(7.4) and (7.5) (with $\Delta$ replaced by $K_\alpha$ and $n$  
replaced by $n_\alpha$).  

For each $\alpha\in\D$ choose $(f_\alpha^n)_{n=1}^\infty    
\subseteq C(K_\alpha)$,  
$f_\alpha^n \ge 0$ and $\Vert f_\alpha^n\Vert =1$, so that $(f_\alpha^n)  
_{n=1}^\infty$ converges pointwise to $F_{n_\alpha}$ and  is equivalent  
to $(s_n)$ with   
$$|f_\alpha^1 (k)| + \sum_{n=1}^\infty | f_\alpha^{n+1} (k) - f_\alpha^n  
(k)| \le n_\alpha \ \hbox{ for all }\  k\in K_\alpha\ .  
\leqno(7.6)$$  
Let $g_n = (c_\alpha f_\alpha^n)_{\alpha\in\D}$. Clearly $g_n\in Y$  
since $\Vert g_n\Vert  \le \sum_{\alpha\in\D} c_\alpha \le1$ by i).  
Furthermore $\ell_1 \not\hookrightarrow \lb (g_n)\rb $ by the following  
lemma and the fact that for all $\alpha$, ${\ell_1\not\hookrightarrow  
\lb f_\alpha^n : n\in\IN\rb}$.  

\proclaim Lemma 7.5.  
For all $\alpha\in \D$, let $Y_\alpha$ be a closed subspace of $C(K_\alpha)$  
which does not contain $\ell_1$. Let   
$$\widetilde Y_\alpha = \bigl\{ (h_\beta)_{\beta\in\D} \in Y :  
h_\alpha \in Y_\alpha \ \hbox{ and }\ h_\beta \equiv 0\ \hbox{ if }\   
\alpha\ne\beta\bigr\}\ .$$  
Let $Z$ be the closed linear span of $\{\widetilde Y_\alpha :\alpha\in\D\}$.  
Then $Z$ does not contain $\ell_1$.\par  

\proof It is shown in \rf{34} that $Y$ does not contain a sequence  
$(h_n)_{n=1}^\infty = ((h_\alpha^n)_{\alpha\in\D})_{n=1}^\infty$  
which is both equivalent to the unit vector basis of $\ell_1$ and has  
the following property: for all $\alpha_0 \in\D$ there exists $m_0\in\IN$  
so that for $m\ge m_0$ and $\alpha\le_L \alpha_0$, $h_\alpha^m \equiv 0$.  

But if $Z$ contains $\ell_1$, then $Y$ must contain such a sequence   
$(h_n)$. This follows easily from the fact that if $(f_n)_{n=1}^\infty$  
is an $\ell_1$-basis in $Z$, then for all $\varep >0$ and $\alpha_0 \in\D$,  
there exists a normalized block basis $(d_n)_{n=1}^\infty = ((d_\alpha^n)_  
{\alpha\in\D})_{n=1}^\infty$ of $(f_n)$ with $\Vert d_{\alpha_0}^n\Vert_  
{C(K_{\alpha_0})} < \varep$ for all $n$.~\qed  
\medskip  

Thus by \rf{36}   
we may pass to a subsequence of $(g_n)$ which is weak  
Cauchy. By relabeling we assume that $(g_n)$ itself is weak Cauchy  
and converges weak* to $F\in Y^{**}$.  

We next verify (7.2). Let $(h_n)$ be a convex block subsequence of  
$(g_n)$. For $k\in \Delta$ and $h= (h_\alpha)_{\alpha\in\D}\in Y$,  
define $\delta_k(h) = \sum_{\alpha\in\gamma_k} \widetilde h_\alpha(k)$  
where $\gamma_k = \{ \alpha\in \D : k\in K_\alpha\}$. Clearly  
$\delta_k$ is a normalized element of $Y^*$. We shall show that  
$$\left\{\eqalign{&  
\hbox{for all $\alpha \in \D$ there exists $k_\alpha\in K_\alpha$  
and $h= (h_\beta)\in Ba \lb (h_n)\rb $}\cr  
&\hbox{such that $\delta_{k_\alpha} (h) > 7/10$  
and $\delta_k (h) < 3/10$ if $k\in\Delta\setminus K_\alpha$.}\cr}  
\right.\leqno(7.7)$$

As in \rf{33} this implies $\lb (h_n)\rb^*$ is nonseparable. Indeed by (7.7)  
we can choose $(h^\alpha)_{\alpha\in\D} \subseteq Ba \lb (h_n)\rb$ and a   
collection of basic clopen sets $(K'_\alpha)_{\alpha\in\D}$ in $\Delta$  
such that for all $\alpha\in\D$,  
\smallskip  
\iitem{a)} $K'_{\alpha,0} \cap K'_{\alpha,1} = \emptyset$,  
\iitem{b)} $K'_{\alpha,\varep} \subseteq K'_\alpha$ for $\varep = 0,1$ and  
\iitem{c)} $\delta_k (h_\alpha) > 7/10$ for $k\in K'_\alpha$ and  
\iitem{} $\delta_k (h_\alpha) < 3/10$ for $k\notin K'_\alpha$.  
\smallskip  

\noindent For each branch (a maximal subset linearly ordered  
by $<$ ) $\gamma$ in  
$\D$ choose $k_\gamma \in\bigcap_{\alpha\in\gamma} K'_\alpha$. By a)  
and b) $k_\gamma$ is well defined and $k_\gamma \ne  k_{\gamma'}$ if  
$\gamma\ne \gamma'$. By c), $\Vert (\delta_{k_\gamma} - \delta_{k_{\gamma'}})  
\big\vert_{\lb (h_n)\rb}\Vert > 2/5$ if $\gamma\ne \gamma'$.  

We return to the proof of (7.7). Fix $\alpha\in\D$    
and set $h_n = (h_\beta^n)_  
{\beta\in\D}$. Since $(h_\alpha^n)_{n=1}^\infty$ is a convex block  
subsequence of $(c_\alpha f_\alpha^n)_{n=1}^\infty$,   
$(h_\alpha^n)_{n=1}^\infty$  
converges pointwise to $c_\alpha F_{n_\alpha}$. Thus by (7.5) and (7.6)  
we may assume (by passing to a subsequence and relabeling, if necessary)  
that there exist $\varep_i = \pm 1$ $(1\le i\le n_\alpha)$ and $k_\alpha  
\in K_\alpha$ such that  
$$n_\alpha \ge \sum_{i=1}^{n_\alpha} c_\alpha^{-1} \varep_i   
(h_\alpha^{i+1} - h_\alpha^i) (k_\alpha) > n_\alpha -1\ .$$  
Let $h= n_\alpha^{-1} c_\alpha^{-1} \sum_{i=1}^{n_\alpha} \varep_i  
(h_{i+1} -h_i) \equiv (h_\beta)_{\beta \in \D}$. Thus  
$1\ge h_\alpha (k_\alpha) > 1-n_\alpha^{-1} > 9/10$ by iv).   

Furthermore by applying (7.6) to each $\beta < \alpha$ we have from ii)  
$$\eqalign{\sum_{\beta <\alpha} \tilde h_\beta (k_\alpha)   
& \le \sum_{\beta <\alpha} n_\alpha^{-1} c_\alpha^{-1} c_\beta n_\beta\cr  
\noalign{\vskip6pt}  
&\le c_\alpha^{-1} n_\alpha^{-1} \sum_{\beta < \alpha} n_\beta < 1/10\ .\cr}$$  
By the triangle inequality and the definition of $h$,  
$$\eqalign{ \sum_{\beta >\alpha} \tilde h_\beta (k_\alpha) 
&\le c_\alpha^{-1} n_\alpha^{-1} \sum_{\beta > \alpha} 2c_\beta n_\alpha\cr  
\noalign{\vskip6pt}  
&= 2c_\alpha^{-1}\sum_{\beta >\alpha} c_\beta   
< 1/10\ \hbox{ (by iii) ).}\cr}$$  
Thus $\delta_{k_\alpha}(h) > 9/10 - 2/10 = 7/10$ which proves the first  
part of (7.7).  

Let $k\in\Delta\setminus K_\alpha$ be fixed. There exists a unique  
$\alpha_0\in \D$ $(\alpha_0\ne \alpha)$ with the same length as   
$\alpha_0$, $|\alpha| = |\alpha_0|$, such  that $k\in K_{\alpha_0}$.  
The calculations above yield $\sum_{\beta <\alpha_0} \tilde h_\beta (k)  
+ \sum_{\beta>\alpha_0} \tilde h_\beta (k) < 2/10$. If $\alpha_0 <_L \alpha$  
then by (7.6)  
$$\eqalign{0\le h_{\alpha_0} (k) & = n_\alpha^{-1} c_\alpha^{-1}  
\sum_{i=1}^{n_\alpha} \varep_i (h_{\alpha_0}^{i+1} - h_{\alpha_0}^i)(k)\cr  
&\le n_\alpha^{-1} c_\alpha^{-1} c_{\alpha_0} n_{\alpha_0} \le 1/10  
\ \hbox{ (by vi) ).}\cr}$$  
If $\alpha <_L \alpha_0$ then we have (from the equality above) that  
$$0 \le h_{\alpha_0} (k) \le n_\alpha^{-1} c_\alpha^{-1} c_{\alpha_0}  
2n_\alpha = 2c_{\alpha_0} c_\alpha^{-1} < 1/10$$  
by v) ). It follows that $\delta_k(h) < 3/10$ which completes the proof  
of (7.7).  

Finally, we verify (7.3). Let $S= \lb \alpha,\beta\rb  \equiv \{\gamma \in\D  
\mid \alpha \le \gamma\le\beta\}$ be a {\it finite segment\/} in $\D$.  
For $k\in K_\beta$ and $f\in Y$ we set $\delta_{S,k}(f) = \sum_{\gamma  
\in S} \tilde f_\gamma (k)$. $\delta_{S,k}(f)$ is defined similarly if  
$S= \lb \alpha,\infty) \equiv \{ \gamma\in \D : \alpha\le\gamma\}$ is an  
infinite segment and $k\in \bigcap_{\beta\in S} K_\beta$. Define  
$$\eqalign{  
K &= \biggl\{ \sum_{i=1}^\infty a_i\delta_{S_i,k_i} : (a_i)_1^\infty \in Ba  
(\ell_2)\ ,\   
	(S_i)_1^\infty \ \hbox{ are disjoint  segments and }\cr   
&\qquad	k_i \in \bigcap_{\beta \in S_i} K_\beta \ \hbox{ for every } \ i  
	\biggr\}\ .\cr}$$  
From the definition of the norm in $Y$ it is clear that $K\subseteq Ba (Y^*)$.  
Furthermore it is easy to check that $K$ is weak* closed and $K$   
1-norms $Y$.  

It remains to show that $F\big|_K \in B_{1/2} (K)$. For $m,n\in\IN$ let  
$g(n,m) \in Y$ be  given by $g(n,m) = (g(n,m)_\beta)_{\beta \in \D}$ where  
$$g(n,m)_\beta = \cases{ g_\beta^n&\quad if $|\beta| \le m$\cr  
\noalign{\vskip4pt}  
0&otherwise.\cr}$$  
Let $y^* = \sum_{i=1}^\infty a_i \delta_{S_i,k_i} \in K$. Then for $m$  
fixed,  
$$\eqalign{\sum_{n=1}^\infty \Big| \big\lb g(n+1,m)-g(n,m)\big\rb (y^*)\Big|    
& = \sum_{n=1}^\infty \Big| \sum_{i=1}^\infty a_i \sum_{\scriptstyle \gamma    
\in S_i\atop\scriptstyle |\gamma| \le m} \big\lb \tilde g_\gamma^{n+1}  
	(k_i) - \tilde g_\gamma^n (k_i)\big\rb\Big| \cr  
\noalign{\vskip6pt}  
&\le \sum_{i=1}^\infty |a_i| \sum_{\scriptstyle \gamma\in S_i\atop  
	\scriptstyle |\gamma| \le m} \sum_{n=1}^\infty  
	| \tilde g_\gamma^{n+1} (k_i) - \tilde g_\gamma^{n+1} (k_i)|\cr  
\noalign{\vskip6pt}  
&\le \sum_{i=1}^\infty |a_i| \sum_{\scriptstyle \gamma\in S_i\atop  
	\scriptstyle|\gamma|\le m} c_\gamma n_\gamma\ \hbox{ (by (7.6) )}\cr  
\noalign{\vskip6pt}  
&\le \sum_{|\gamma| \le m} c_\gamma n_\gamma < \infty\ .\cr}$$  
In particular $(g(n,m))_{n=1}^\infty$ converges pointwise on $K$ to a  
function $G_m \in DBSC (K)$.   

All that remains is to show that $\Vert G_m - F\big|_K \Vert_{C(K)} \to0$  
as $m\to \infty$.  Let $m\in\IN$ be fixed and let $y^* = \sum_{i=1}^\infty a_i  
\delta_{S_i,k_i}\in K$. Then  
$$\eqalign{|G_m (y^*) - F(y^*) |   
&= \Big| \sum_{i=1}^\infty a_i \sum_{\scriptstyle \gamma\in S_i\atop  
	\scriptstyle |\gamma| >m} c_\gamma \tilde F_{n_\gamma} (k_i)\Big|\cr  
\noalign{\vskip6pt}  
&\le \sum_{i=1}^\infty | a_i| \biggl( \sum_{\scriptstyle \gamma\in S_i  
	\atop\scriptstyle |\gamma| >m} c_\gamma\biggr)\ .\cr}$$  

For each $i$ set   
$$b_i^m = \sum_{\scriptstyle \gamma\in S_i\atop\scriptstyle |\gamma| >m}  
c_\gamma\ .$$  
Thus  
$$|G_m(y^*) - F(y^*) | \le \biggl( \sum_{i=1}^\infty (b_i^m)^2\biggr)^{1/2}$$  
by H\"older's inequality. The latter goes to 0 as $m\to\infty$  
by vii).~\qed  

\beginsection{8. Problems.}  

We have previously raised two problems concerning $B_{1/4} (K)$.  

\demo Problem 8.1.  
Let $F\in B_1 (K)$ and $C<\infty$ be such that if $(f_n) \subseteq C(K)$  
is a bounded sequence converging pointwise to $F$, then there exists  
$(g_n)$, a convex block subsequence of $(f_n)$, with spreading model  
$C$-equivalent to the summing basis. Is $F\in B_{1/4}(K)$?  
  
\demo Problem 8.2.  
Let $F\in B_1(K)$ and assume there exists a $C<\infty$ such that if   
$(\varep_i)\subseteq \IR^+$ and  $K_n  
(F,(\varep_i))\ne \emptyset$, then $\sum_1^n \varep_i\le C$.  
Is $F\in B_{1/4} (K)$?  

These problems lead naturally to the following definitions.  
Let $F\in B_1 (K)$.   
$$\eqalign{|F|_I & = \max \left\{ \sup \biggl\{ \sum_{i=1}^m \delta_i  
	: K_m \bigl( F,(\delta_i)\bigr) \ne \emptyset\biggr\}\ ,\   
	\Vert F\Vert_\infty\right\}\ .\cr  
\noalign{\vskip6pt}  
|F|_{I'} & = \max \Bigl\{ \sup \bigl\{ m\delta : K_m (F,\delta) \ne  
	\emptyset \bigr\}\ ,\ \Vert F\Vert_\infty\Bigr\}\ .\cr  
\noalign{\vskip6pt}  
|F|_S & = \inf \biggl\{ C : \ \hbox{ there exist $(f_n) \subseteq C(K)$  
	converging pointwise to $F$}\cr  
&\qquad \hbox{with for all }\ (a_i)_1^k\subseteq \IR\ ,\   
	\lim_{\scriptstyle n_1\to\infty\atop\scriptstyle n_1<\cdots < n_k}  
	\Big\Vert \sum_{i=1} ^k a_if_{n_i}\Big\Vert   
	\le C\Big\Vert \sum_1^k a_is_i\Big\Vert \biggr\}\ .\cr}$$  

\demo Remark 8.3.  
We do not know if $|F|_I$ or $|F|_{I'}$ are norms. It is clear that $|F|_S$  
is a norm and also that  
$$\Vert F\Vert_\infty \le |F|_{I'} \le |F|_I \le |F|_S \le |F|_{1/4}  
\le |F|_D$$  
($|F|_S \le |F|_{1/4}$ follows from the proof of Theorem B.)  Furthermore,  
using the series criterion for completeness, it is easy to show that  
$(\{ F\in B_1 (K): |F|_S <\infty\},\, |\cdot|_S)$ is a Banach space.  

\demo Problem 8.4.  
Are $|\cdot|_I$ and $|\cdot |_S$ equivalent? Are $|\cdot|_S$ and  
$|\cdot |_{1/4}$ equivalent?  

The  solution  of Problem 8.4 would of course solve Problems~8.1 and 8.2.  
Furthermore an affirmative answer to Problem~8.2 would yield an affirmative  
answer to Problems~8.1 and 8.4.  

\proclaim Proposition 8.5.  
$|\cdot|_I$ and $|\cdot|_{I'}$ are not (in general) equivalent. \par  

\proof Define  $F: \lb 0,1\rb^\omega\to \IR$ as follows:  
\smallskip  
\item{} If $t_0\ne 0$ let  
$$F(t_0,t_1,\ldots) = \sin t_0^{-1}\ .$$  
\item{} If $t_0= t_1 = \cdots = t_r = 0\ne t_{r+1}$, set  
$$F(t_0,t_1,\ldots) = {1\over r+2} \sin t_r^{-1}\ .$$  
\smallskip  

\noindent It's easy to see that $\osc (F; (0,t_1,t_2,\ldots)) =2$ for all  
$t_1,t_2,\ldots \in \lb 0,1\rb$ and so  
$$K_1 (F,\varep) = \{ 0\} \times \lb 0,1\rb^{\omega\setminus \{ 0\}}$$  
whenever $0<\delta <2$. Similar calculations show that if $r=\lbrack\!\lbrack  
{2\over\varep}\rbrack\!\rbrack$ then  
$$K_r (F,\varep) = \{ 0\}^r \times \lb 0,1\rb^{\omega\setminus r}$$  
and $K_{r+1} (F,\varep) = \emptyset$. Thus $K_m (F,\varep) \ne\emptyset$  
implies $m\varep \le2$. On the other hand, for $m\ge1$,  
$$K_m \left( F, \Bigl( 2,1,{2\over3},\cdots, {2\over m}\Bigr)\right) =  
\{ 0\}^m \times \lb 0,1\rb^{\omega\setminus m}\ .\eqno\blackbox$$  

We conclude by mentioning some further problems for study, some of which have  
been raised above.  

\demo Problem 8.6.  
Classify (or give useful sufficient conditions) for a function $F\in B_1(K)$  
to govern $\{ X: X^*$ is separable  and $\dim X = \infty\}$. In particular  
is $F\in B_{1/4} (K)\setminus C(K)$ a sufficient condition?  

\demo Problem 8.7.  
Classify those $F\in B_1(K)$ which govern $\{ \ell_1\}$, which govern  
$\{ c_0\}$, which govern $\{ X: X$ is reflexive$\}$ or which govern  
$\{ X: X$ is quasi-reflexive$\}$.   

We note that if $X$ is a Polish Banach space ({\it i.e.}, $Ba(X)$ is  
Polish in the weak topology) then Edgar and Wheeler \rf{14} have shown  
that $X$ is hereditarily reflexive (see also \rf{37} and \rf{18}).  
Bellenot \rf{5} and Finet \rf{15} have independently extended this result  
by showing that whenever $X$ is Polish, if $x^{**} \in X^{**}\setminus  
X$ then $x^{**}|_{Ba (X^*)}$ strictly governs the class of quasi-reflexive  
spaces of order~1.   
\bigskip  

\baselineskip=12pt  
\frenchspacing  
\centerline{\bf References}  

\item{1.} A. Andrew, {\it Spreading basic sequences and subspaces of  
James' quasi-reflexive space}, Math. Scand. {\bf48} (1981), 109--118.   

\item{2.} P. Azimi and J.N. Hagler,  
{\it Examples of hereditarily $\ell^1$ Banach spaces failing the Schur  
property}, Pacific J. Math. {\bf122} (1986), 287--297.  

\item{3.} R. Baire, {\it Sur les Fonctions des Variables R\'eelles},  
Ann. di  Mat. {\bf 3} (1899), 1--123.  

\item{4.} B. Beauzamy and J.-T. Laprest\'e,  
{\it Mod\`eles \'etal\'es des espaces de Banach},   
Travaux en Cours, Hermann, Paris (1984).  

\item{5.} S. Bellenot,  
{\it More quasi-reflexive subspaces}, Proc. AMS {\bf 101} (1987),  
693--696.  

\item{6.} S. Bellenot, R. Haydon and E. Odell,  
{\it Quasi-reflexive and tree spaces constructed in the spirit of R.C.~James},  
Contemporary Math. {\bf85} (1989), 19--43.  

\item{7.} C. Bessaga and A. Pe{\l}czy\'nski,  
{\it On bases and unconditional convergence of series in Banach spaces},  
Stud. Math. {\bf17} (1958), 151--164.  

\item{8.} J. Bourgain,  
{\it On convergent sequences of continuous functions},  
Bull. Soc. Math. Bel. {\bf32} (1980), 235--249.  

\item{9.} J. Bourgain,  
{\it Remarks on the double dual of a Banach space},  
Bull. Soc. Math. Bel. {\bf32} (1980), 171--178.  

\item{10.} J. Bourgain, unpublished notes.  

\item{11.} P.G. Casazza and T.J. Shura,  
{\it Tsirelson's Space}, Springer-Verlag Lecture Notes in Mathematics,  
{\bf 1363} (1989).   

\item{12.} W.J. Davis, T. Figiel, W.B. Johnson and A. Pe{\l}czy\'nski,  
{\it Factoring weakly compact operators},  
J. Funct. Anal. {\bf17} (1974), 311--327.  

\item{13.} J. Elton,  
{\it Extremely weakly unconditionally convergent series},  
Israel J. Math. {\bf40} (1981), 255--258.  

\item{14.} G.A. Edgar and R.F. Wheeler,  
{\it Topological properties of Banach spaces},  
Pacific J. Math. {\bf115} (1984), 317--350.  

\item{15.} C. Finet,  
{\it Subspaces of Asplund Banach spaces with the point of continuity  
property}, Israel J. Math. {\bf 60} (1987), 191--198.  

\item{16.} V. Fonf,  
{\it One property of Lindenstrauss-Phelps spaces},  
Funct. Anal. Appl. (English trans.) {\bf13} (1979), 66--67.  

\item{17.} T. Figiel and W.B. Johnson,  
{\it A uniformly convex Banach space which contains no $\ell_p$},  
Comp. Math. {\bf29} (1974), 179--190.  

\item{18.} N. Ghoussoub and B. Maurey,  
{\it $G_\delta$-embeddings in Hilbert space},  
J. Funct. Anal. {\bf61} (1985), 72--97.  
 
\item{19.} \refrule \ ,  
{\it $G_\delta$-embeddings in Hilbert space II},  
J. Funct. Anal. {\bf 78} (1988), 271--305.  
 
\item{20.} \refrule\ ,  
{\it $H_\delta$-embeddings in Hilbert space and optimization   
on $G_\delta$ sets},  
Memoirs Amer. Math. Soc.  {\bf62} (1986), number 349.  
 
\item{21.} \refrule\ ,  
{\it A non-linear method for constructing certain basic sequences   
in Banach spaces}, preprint.  

\item{22.} N. Ghoussoub, G. Godefroy, B. Maurey and W. Schachermayer,  
{\it Some topological and geometrical structures in Banach spaces},  
preprint.  

\item{23.} F. Hausdorff,  
``Set Theory'', Chelsea, New York (1962).  

\item{24.} R. Haydon and B. Maurey, {\it On Banach spaces with strongly   
separable  types}, J.~London Math. Soc. {\bf33} (1986), 484--498.  

\item{25.} A.S. Kechris and A. Louveau,  
{\it A classification of Baire class 1 functions}, preprint.  

\item{26.} J.L. Krivine and B. Maurey, {\it Espaces de Banach stables},  
Israel J. Math. {\bf39} (1981), 273--295.  

\item{27.} J. Lindenstrauss and L. Tzafriri,  
``Classical Banach spaces'',  
Springer-Verlag Lecture Notes in Math. {\bf 338}, Berlin (1973).  

\item{28.} \refrule\ ,  
``Classical Banach spaces II'',  
Springer-Verlag, Berlin (1977).  

\item{29.} A.A. Milutin,  
{\it Isomorphisms of spaces of continuous functions on compacta of   
power continuum},  
Tieoria Funct. (1966), 150--166 (Russian).  

\item{30.} S. Mazurkiewicz and W. Sierpinski,  
{\it Contribution \`a la topologie des ensembles d\'e nombrales},  
Fund. Math. {\bf1} (1920), 17--27.  

\item{31.} A. Pe{\l}czy\'nski,  
{\it A note on the paper of I.~Singer ``Basic sequences and reflexivity of  
Banach spaces''},  
Studia Math. {\bf21} (1962), 371--374.  

\item{32.} E. Odell,   
{\it A nonseparable Banach space not containing   
a subsymmetric basic sequence},  
Israel J. Math. {\bf52} (1985), 97--109.  

\item{33.} \refrule\ ,  
{\it Remarks on the separable dual problem},  
Proceedings of Research Workshop on Banach Space Theory (ed. by B.-L.Lin),  
The University of Iowa (1981), 129--138.  

\item{34.} \refrule\ ,  
{\it A normalized weakly null sequence with no  shrinking subsequence in  
a Banach space not containing $\ell_1$},  
Comp. Math. {\bf41} (1980), 287--295.  

\item{35.} E. Odell and H. Rosenthal,  
{\it A double-dual characterization of separable Banach spaces not  
containing $\ell_1$},  
Israel J. Math. {\bf20} (1975), 375--384.  

\item{36.} H. Rosenthal,  
{\it A characterization of Banach spaces containing $\ell_1$},  
Proc.  Nat. Acad. Sci. U.S.A. {\bf71} (1974), 2411--2413.  

\item{37.} \refrule\ ,  
{\it Weak*-Polish Banach spaces},  
J. Funct. Anal. {\bf76} (1988), 267--316.  

\item{38.} \refrule\ ,  
{\it Some remarks concerning unconditional basic sequences},  
Longhorn Notes, University of Texas, (1982-83), 15--48.  

\item{39.} A. Sersouri,  
{\it A note on the Lavrientiev index for the quasi-reflexive Banach spaces},  
Contemporary Math. {\bf 85} (1989), 497--508.

\bigskip  
\bigskip  
{\vbox{\halign{\hbox{#}\hfil\qquad&\hbox{#}\hfil\qquad&\hbox{#}\hfil\cr  
R. Haydon&E. Odell&H. Rosenthal\cr  
Brasenose College&The University of Texas at Austin&The University of Texas  
	at Austin\cr  
Oxford OX1 4AJ&Austin, Texas 78712&Austin, Texas 78712\cr  
England&U.S.A.&U.S.A.\cr}}}

\bigskip  
\rightline{September 6, 1990}  
\end